\documentclass[12pt]{amsart}
\usepackage{color}
\usepackage[T1]{fontenc}
\usepackage[utf8]{inputenc}
\usepackage{hyperref,amssymb,url,upref,verbatim,xspace,mathrsfs}
\RequirePackage[all,ps,cmtip]{xy}
\newdir^{ (}{{}*!/-5pt/\dir^{(}}
\newdir_{ (}{{}*!/-5pt/\dir_{(}}
\RequirePackage[mathscr]{eucal}
\usepackage{mathrsfs}\let\mathcal\mathscr
\theoremstyle{plain}

\newtheorem{theorem}{Theorem}[section]
\newtheorem{proposition}[theorem]{Proposition}
\newtheorem{lemma}[theorem]{Lemma}
\newtheorem{corollary}[theorem]{Corollary}
\newtheorem{assumption}[theorem]{Assumption}
\theoremstyle{definition}
\newtheorem{definition}[theorem]{Definition}
\newtheorem{remark}[theorem]{Remark}
\newtheorem{example}[theorem]{Example}
\numberwithin{equation}{section}

\newcounter{toto}
\def\thetoto{\arabic{toto}}
\let\oldmarginpar\marginpar
\def\marginpar#1{\refstepcounter{toto}\textsuperscript{\textup{[\thetoto]}}\oldmarginpar{\footnotesize\textsuperscript{[\thetoto]}\,#1}}

\def\mainmatter{\renewcommand{\baselinestretch}{1.1}\normalfont}

\makeatletter
\@namedef{subjclassname@2020}{%
  \textup{2020} Mathematics Subject Classification}
\makeatother

\newcommand{\CF}{{\mathrm{CF}}}
\newcommand{\C}{\mathbb{C}}

\newcommand{\R}{\mathbb{R}}

\newcommand{\Z}{\mathbb{Z}}

\newcommand{\bD}{\boldsymbol{D}}
\newcommand{\shhom}{\mathcal{H}\!\mathit{om}}
\DeclareMathOperator{\rh}{\mathit{R}\shhom}\let\Rhom\rh
\DeclareMathOperator{\tho}{\mathit{T}\shhom}

\DeclareMathOperator{\RH}{RH}

\newcommand{\rb}{\mathrm{b}}

\newcommand{\alg}{\mathrm{alg}}

\newcommand{\coh}{\mathrm{coh}}
\newcommand{\hol}{\mathrm{hol}}
\newcommand{\rhol}{\mathrm{rhol}}

\newcommand{\Mod}{\mathrm{Mod}}

\newcommand{\Pic}{\mathrm{Pic}}

\newcommand{\cc}{{\C\textup{-c}}}
\newcommand{\rc}{{\R\textup{-c}}}
\newcommand{\good}{\textup{-}\rm{good}}

\newcommand{\wrc}{{\textup{w-}\R\textup{-c}}}
\newcommand{\XS}{X\times S}

\newcommand{\DXS}{\shd_{\XS/S}}

\DeclareMathOperator{\rD}{\mathsf{D}}
\DeclareMathOperator{\rK}{\mathsf{K}}
\DeclareMathOperator{\DR}{DR}

\DeclareMathOperator{\pDR}{{}^\mathrm{p}DR}

\DeclareMathOperator{\Hom}{Hom}
\DeclareMathOperator{\id}{Id}\let\Id\id

\DeclareMathOperator{\Sol}{Sol}
\DeclareMathOperator{\pSol}{{}^\mathrm{p}Sol}
\DeclareMathOperator{\supp}{Supp}

\let\ov\overline
\let\epsilon\varepsilon

\let\setminus\smallsetminus

\def\loccit{loc.\kern3pt cit.{}\xspace}
\def\cf{cf.\kern.3em}
\def\eg{e.g.\kern.3em}

\def\resp{\text{resp.}\kern.3em}

\newcommand{\Dq}{{}_{\scriptscriptstyle\mathrm{D}}q}

\newcommand{\cbbullet}{{\raisebox{1pt}{$\sbullet$}}}
\newcommand{\sbullet}{{\scriptscriptstyle\bullet}}
\newcommand{\pOS}{p^{-1}_X\sho_S}
\newcommand{\puOS}{p^{-1}_U\sho_S}

\def\sha{\mathcal{A}}

\def\shd{\mathcal{D}}
\def\she{\mathcal{E}}
\def\shf{\mathcal{F}}\let\cF F
\def\shg{\mathcal{G}}\let\cG G
\def\shh{\mathcal{H}}
\def\shi{\mathcal{I}}

\def\shk{\mathcal{K}}
\def\shl{\mathcal{L}}
\def\shm{\mathcal{M}}

\def\sho{\mathcal{O}}
\def\shp{\mathcal{P}}

\def\sht{\mathcal{T}}

\newcommand{\RedefinitSymbole}[1]{
\expandafter\let\csname old\string#1\endcsname=#1
\let#1=\relax
\newcommand{#1}{\csname old\string#1\endcsname\,}
}
\RedefinitSymbole{\forall} \RedefinitSymbole{\exists}

\def\to{\mathchoice{\longrightarrow}{\rightarrow}{\rightarrow}{\rightarrow}}

\def\To#1{\mathchoice{\xrightarrow{\textstyle\kern4pt#1\kern3pt}}{\stackrel{#1}{\longrightarrow}}{}{}}

\def\isom{\stackrel{\sim}{\longrightarrow}}

\begin{document}

\baselineskip=17pt

\title{On relative constructible sheaves and integral transforms}

\author[L. Fiorot]{Luisa Fiorot}
\address[L. Fiorot]{Universit\`a degli Studi di Padova\\
Dipartimento di Matematica ``Tullio Levi-Civita'' 
Via Trieste, 63\\
35121 Padova Italy}
\email{luisa.fiorot@unipd.it}

\author[T. Monteiro Fernandes]{Teresa Monteiro Fernandes}
\address[T. Monteiro Fernandes]{Centro de Matem\'atica, Aplica\c{c}\~{o}es Funda\-men\-tais e Investiga\c c\~ao Operacional and Departamento de Matem\' atica da Faculdade de Ci\^en\-cias da Universidade de Lisboa, Bloco C6, Piso 2, Campo Grande, 1749-016, Lisboa
Portugal}
\email{mtfernandes@fc.ul.pt}

\date{}

\begin{abstract}
The aim of this note is threefold. 
The first is to 
obtain a simple characterization of relative constructible sheaves when the parameter space is projective. The second is to study the relative Fourier-Mukai  for relative constructible sheaves and  for relative regular holonomic $\shd$-modules and prove they induce relative equivalences of categories.
 The third is to introduce and study the notions of relative constructible functions and relative Euler-Poincar\'e index. We prove that the relative Euler-Poincar\'e index provides an isomorphism between the Grothendieck group of the derived category of complexes with bounded relative $\mathbb R$-constructible cohomology and the ring of relative constructible functions.
\end{abstract}

\maketitle
\tableofcontents
\mainmatter

\section*{Introduction}

Let $X$ be  a real analytical manifold of dimension $d_X$, let $S$  be a complex manifold and let 
$p_X: X\times S\to S$ be  the projection. Let $\coh(\sho_S)$ (resp. $\rD^\rb_{\coh}(\sho_S)$) denote
the category of $\sho_S$-coherent  modules (resp. the triangulated category of complexes of $\sho_S$-modules with bounded $\sho_S$-coherent cohomology). Motivated by the study of the relative Riemann-Hilbert correspondence, the notions of relative $\R\textup{-}$constructibility as well as that of relative 
$\C\textup{-}$constructibility, when $X$ is a complex manifold, were introduced in \cite{MFCS1}. Let $\Mod_{\rc}(\pOS)$ (resp. $\Mod_{\cc}(\pOS)$) be the category of $S\textup{-}\R\textup{-}$(resp. $S\textup{-}\C\textup{-}$)constructible sheaves of $\pOS$-modules.

Let $\rD^\rb_{\rc}(\pOS)$ (resp. $\rD^\rb_{\cc}(\pOS)$) denote the triangulated category of complexes of $\pOS$-modules with bounded $S\textup{-}\R\textup{-}$constructible cohomologies (resp. the triangulated category of complexes of $\pOS$-modules with bounded $S\textup{-}\C\textup{-}$constructible cohomologies in the complex case)  introduced in \cite{MFCS1}.

As proved in \cite[Cor.\,A.4]{FMFS2}, the category $\rD^\rb_{\rc}(\pOS)$ is equivalent to the category whose objects are the weakly $S$-$\R$-constructible complexes with bounded $S$-$\R$-constructible cohomology. 

The space $S$ is considered as a space of parameters, so that, when $S=\{pt\}$ (the absolute case), one recovers the category $\rD^\rb_{\rc}(X)$ whose objects are complexes of sheaves of $\mathbb C$-vector spaces with bounded $\R$-constructible cohomology.  Let $\rD^\rb(\Mod_{\rc}(X))$ denote the bounded derived category of the abelian category $\Mod_{\rc}(X)$. Kashiwara and Schapira proved in \cite[Th.\,8.1.11]{KS1} that the canonical functor 
\begin{equation}\label{EN2}\rD^\rb(\Mod_{\rc}(X))\to \rD^\rb_{\rc}(X)
\end{equation} is an equivalence of categories.

When $X=\{pt\}$ one recovers the category $\rD^\rb_{\coh}(\sho_S)$ whose objects are  complexes of $\sho_S$-modules with bounded coherent cohomology. 
Let $\rD^\rb(\coh(\sho_S))$ denote the bounded derived category of the abelian category ${\coh}(\sho_S)$. 
The canonical functor 
\begin{equation}\label{E200} E:\rD^\rb(\coh(\sho_S))\to \rD^\rb_{\coh}(\sho_S)
\end{equation} is not in general known to be an equivalence so that the functor 
\begin{equation}\label{EN1}
 \rD^\rb(\Mod_{\rc}(\pOS))\to\rD^\rb_{\rc}(\pOS)
\end{equation}
 is not in general known to be an equivalence, contrary to the absolute case.
The first main result in this work is Theorem~\ref{Cor1} where we prove that \eqref{EN1} is an equivalence under the assumption that $S$ is a projective algebraic complex manifold. Conversely,  in Proposition \ref{PN3} we prove that a necessary condition (on $S$) for \eqref{EN1} being an equivalence for arbitrary $X$ is precisely that \eqref{E200} is an equivalence. This last condition is very natural in the sense that any projective algebraic complex manifold satisfies it as proved in Proposition \ref{Remalgan}.
Having proved in \ref{pre1} the necessary functorial properties of relative $\R$-constructible complexes and of the relative Riemann-Hilbert functor (stability under pullback and pushforward on the parameter space), in Section \ref{IT} we embark on the study of relative integral transforms for relative constructible sheaves and for relative $\shd$-modules. An example of the latter (\cf Example \ref{ERIT}) is given by the generalized Fourier-Mukai transforms. Let us consider an  abelian variety $A$ and the associated  generalized Fourier-Mukai integral transform $\rm{FM}$  recalled in \cite[Intro. (iii)]{FMFS2}, due to Laumon \cite{L} and Rothstein \cite{Rot96}. For each object $\shm$ in $\rD^\rb_{\coh}(\shd_A)$, $\rm{FM}(\shm)$ is an object of $\rD^\rb_{\coh}(\sho_{A^{\sharp}})$,  where $A^{\sharp}$ is the  the moduli space of line bundles with integrable connection on A (cf. \cite{L} and \cite{Sch15} for details) which is a quasi-projective algebraic variety equipped with a natural structure of abelian group. The kernel $\mathcal{P}$ of $\rm{FM}$ is a flat relative connection, thus a relative regular holonomic $\shd_{A\times A^{\sharp}/A^{\sharp}}-$ module.

We then pursue our study by considering the relative transformations arising from the Fourier-Mukai transform in the parameter space. Denoting by $\hat{A}$ the dual abelian variety, it is easy to show that  the relative Fourier-Mukai transform provides equivalences of categories $\rD^\rb_{\rc}(p^{-1}_X\sho_A)\to \rD^\rb_{\rc}(p^{-1}_X\sho_{\hat{A}})$ (respectively $\rD^\rb_{\cc}(p^{-1}_X\sho_A)\to \rD^\rb_{\cc}(p^{-1}_X\sho_{\hat{A}})$). 
On the other hand, having fixed the parameter manifold $S$, applying the results in (\cite{D-S}), in the case $X=\mathbb{P}$, the complex projective space of dimension $n$, and if $\mathbb{P}^*$ is its dual, we show that the relative projective duality induces an equivalence of categories $\rD^\rb_{\cc}(p^{-1}_{\mathbb{P}}\sho_S)\to \rD^\rb_{\cc}(p^{-1}_{\mathbb{P}^*}\sho_{S})$.
As a consequence, by the relative Riemann-Hilbert correspondence (\cite{FMFS2}), we conclude that if $A$ is a complex abelian variety then the bounded  derived categories of complexes with relative regular holonomic cohomologies respectively on $\mathbb{P}\times A$ and $\mathbb{P}^*\times \hat{A}$ are equivalent. 
Moreover the relative Riemann-Hilbert correspondence allows to translate all these equivalences in terms of integral transforms for relative $D$-modules, determining the respective kernels.

In the third part of this note (Section \ref{SCF}, independent of the first ones), our motivation was the conceptual need of a notion of relative constructible functions.
We recall that, to the notion of $\R$-constructibility in the absolute case, one associates that of constructible function (taking values in $\Z$) which has been treated in several works (for our main purposes we refer to \cite{KS1}, \cite{S1} and \cite{S2}) with striking applications to  topological data analysis (TDA) (\cite{V}). One of the key results of this theory is  the isomorphism between the Grothendieck ring $\rK_{\rc}(X)$ of the category $\rD^\rb_{\rc}(X)$ and the ring $\CF(X)$ of constructible functions on $X$ via the Euler-Poincar\'e index $\chi$ (\cite[Th.\,9.7.1]{KS1}).
 It became clear that $\Z$ was to be regarded as the Grothendieck ring $\rK_0(\mathbb C\textup{-}{\mathrm{Vect}})$ of the abelian category of finitely generated $\C$-vector spaces. Then, a natural candidate to play the role of $\mathbb Z$ in the relative setting is the Grothendieck group (and unital ring) $\rK_0(S)$ of the triangulated category $\rD^\rb_{\coh}(\sho_S)$.

After stating a natural definition of relative constructible function, in Theorem~\ref{T1fSc} we obtain
  the isomorphism $\chi$ in  the relative setting. Let us be more precise:

From one side the natural generalization of $\rK_{\rc}(X)$ is the Grothendieck group $\rK_0(\rD^\rb_{\rc}(\pOS))$  of
$\rD^\rb_{\rc}(\pOS)$. 
On the other side one naturally adapts the notion of constructible function introducing the $S$-constructible functions on $X$ (Definition~\ref{DfSc}) 
(\S~\ref{S}) by requiring similar conditions to \cite [(9.7.1)]{KS1} but with $\Z$ replaced by $\rK_0(S)$. We then introduce the ring $\CF^S(X)$ as the ring of global $S$-constructible functions.

We recall the following key fact: for each $x\in X$ and each $F\in\rD ^\rb_{\rc}(\pOS)$, $F|_{\{x\}\times S}$ is an object of $\rD^\rb_{\coh}(\sho_S)$. 

Then the natural generalization of the local Euler-Poincar\'e index $\chi$ is the relative Euler-Poincar\'e index  $\chi^S$
(Definition~\ref{Def:EPI}) obtained by assigning to $F$ the class of $F|_{\{x\}\times S}$ in $\rK_0(S)$.

Note that, in the case of $X=\{pt\}$, $\chi^S$   reduces to 
the identity.

The proof of Theorem~\ref{T1fSc} is stepwise similar to that of \cite[Th. 9.7.1]{KS1}.

A particular interesting case is when $F$ is the holomorphic solution  complex $\Sol \shm$ of a bounded complex $\shm$ of $\DXS$-modules with holonomic cohomologies (i.e $\shm$ is an object of $\rD^\rb_{\hol}(\DXS)$, see Ex.\ref{EX22}). $F$ is then $S$-$\C$-constructible and we thereby obtain a relative local index $\chi^S$ on $\rD^\rb_{\hol}(\DXS)$ setting $$\chi^S(\shm):=\chi^S(\Sol \shm)$$ which coincides with Kashiwara's local index in the absolute case.
The existence of  a relative variant of Kashiwara's local index theorem (\cf for instance \cite[Th. 6.3.1]{K8}) for holonomic $\shd_X$-modules remains conjectural however.

 We warmly thank Pierre Schapira for corrections and suggestions on a previous version. We are grateful to Haohao Liu for pointing us a gap. We also thank Amnon Neeman, Andrea D'Agnolo, Claude Sabbah and Thomas Kr\"amer for useful explanations. Finally we are grateful to the referee for her/his  questions, corrections and criticism.
\section{Relative constructible sheaves} \label{pre}
\subsection{Background on relative constructible sheaves and the relative Riemann-Hilbert correspondence}\label{pre1}Let  $X$ be a real analytic manifold and let $S$ be a complex manifold. Following \cite[\S 2]{MFCS1}, a sheaf $F$ of $\pOS$-module  is called
$S$-locally constant coherent if, for each point $(x_0,s_0)\in X\times S$, there exists a neighborhood $U=V_{x_0}\times T_{s_0}$ and
a coherent sheaf $G^{(x_0,s_0)}$ of $\sho_{T_{s_0}}$-modules such that $F|_{U}\cong
p_{V_{x_0}}^{-1}(G^{(x_0,s_0)})$. 
As explained in \cite[\S.~2.2]{MFCS1}, if we choose a sufficiently small contractible neighborhood $V_{x_0}$ of $x_0$ then $U=V_{x_0}\times S$ and $G^{(x_0,s_0)}=F|_{\{x_0\}\times S}$ is a coherent $\sho_S$-module.

We define $\rD^\rb_{{\mathrm {lc}}\; {\mathrm {coh}}}(\pOS)$ to be the full
subcategory of $\rD^\rb(\pOS)$ whose complexes have $S$-locally
constant coherent cohomologies (hence $F|_{\{x_0\}\times S}\in \rD^\rb_{\mathrm {coh}}(\sho_S)$).

We denote by $\Mod_\wrc(\pOS)$ (resp. $\Mod_\rc(\pOS)$) the full subcategory of $\Mod(\pOS)$ whose 
objects $F$ 
are  such that, for a suitable $\mu$-stratification
$(X_\alpha)$ of $X$, 
$F|_{X_{\alpha}\times S}$ is $S$-locally constant
(resp. $S\textup{-}$locally constant coherent).
We call these objects
weakly  $S\textup{-}\R\textup{-}$constructible sheaves (resp. $S\textup{-}\R\textup{-}$constructible sheaves).
When $X$ is a complex analytic manifold and the stratification $(X_{\alpha})$ above is $\C$-analytic, we say that $F$ is $S$-$\C$-constructible. 
Let us denote by $\rD^\rb_{\rc}(\pOS)$ (resp.  $\rD^\rb_\wrc(\pOS)$, resp. $\rD^\rb_{\cc}(\pOS)$) the full
subcategory of $\rD^\rb(\pOS)$ whose objects~$F$ have  
 $S$-$\R$-constructible cohomologies (resp. weakly  $S$-$\R$-constructible cohomologies, resp. $S$-$\C$-constructible coghomologies). 

Let us recall some results on $S$-constructible sheaves:
\begin{itemize}
\item According to \cite[Cor.\,A.4]{FMFS2} the canonical functors induce equivalences of triangulated categories
$\rD^\rb(\Mod_\wrc(\pOS))\simeq \rD^\rb_\wrc(\pOS)$ and
$\rD^\rb_\rc(\Mod_\wrc(\pOS)) \simeq \rD^\rb_\rc(\pOS).$ In order to keep the notation clean in the sequel, we will write $\rD^\rb_\rc(\pOS)$ instead of $\rD^\rb_\rc(\Mod_\wrc(\pOS))$ i.e. the complexes
in $\rD^\rb_\rc(\pOS)$ are always assumed to have weakly $S$-$\R$-constructible entries.

 \item
 Following \cite[Prop.\,2.23]{MFCS1}, 
 the contravariant functor
$\bD_X(\cbbullet)=\rh_{\pOS}(\cbbullet, \pOS)[d_X]$
is a duality in $\rD^\rb_{\rc}(\pOS)$  which induces a duality in $\rD^\rb_{\cc}(\pOS)$. 
\item Following \cite[Prop.\,2.17]{MFCS1},
 let  $F$ in $\rD^\rb_\rc(p_Y^{-1}\sho_S)$ and let $f:Y\to X$  be a morphism of real analytic manifolds 
 which is proper on $\supp_Y(F)$.
Then $Rf_*F$ belongs to  $\rD^\rb_\rc(p_X^{-1}\sho_S)$ where we keep the notation $f$ for the induced morphism $f \times \Id : X\times S\to Y\times S$.
\end{itemize}

For the sake of completeness we include functorial properties which were tacitly assumed but not written in \cite{MFCS1}.
Let $\pi: T\to S$ be a morphism of complex manifolds and let us keep the notation $\pi$ for the induced morphism $\Id\times \pi: X\times T\to X\times S$. 
\begin{proposition}\label{PN1}
For each $F\in\rD^\rb_\rc(p_X^{-1}\sho_S)$ (resp. $F\in\rD^\rb_\cc(p_X^{-1}\sho_S)$), $L\pi^*F:=p_X^{-1}\sho_T\otimes^L_{p_X^{-1}\pi^{-1}\sho_S}\pi^{-1}F$ is a complex in $\rD^\rb_\rc(p_X^{-1}\sho_T)$ (resp. in $\rD^\rb_\cc(p_X^{-1}\sho_T)$) .
\end{proposition}
\begin{proof}
The argument is similar in both cases. We treat the $\R$-constructible case. According to the triangulation Theorem (\cite[Th.~8.2.5]{KS1}) and \cite[Th.~8.3.20]{KS1}, we can choose a $\mu$-stratification $X=\bigsqcup_{\alpha} X_{\alpha}$ such that $F|_{X_{\alpha}\times S}$ is isomorphic to $p_{X_{\alpha}}^{-1}G_{\alpha}$ for some $G_{\alpha}\in\rD^\rb_{\coh}(\sho_S)$. Thus, $L\pi^*F|_{X_{\alpha}\times T}\simeq p_{X_{\alpha}}^{-1}L\pi^*G_{\alpha}$, where $L\pi^*G_{\alpha}\in\rD^\rb_{\coh}(\sho_T)$.
\end{proof}
\begin{proposition}\label{PN2}
Let us be given $G\in\rD^\rb_\rc(p_X^{-1}\sho_T)$ (resp. $G\in \rD^\rb_\cc(p_X^{-1}\sho_T)$) and let us assume that $\pi$ is proper on  $\supp_TG:=\ov{p_X(\supp G)}$. Then, $R\pi_*G$ is a complex 
in $\rD^\rb_{\rc}(p_X^{-1}\sho_S)$ (resp. in $\rD^\rb_{\cc}(p_X^{-1}\sho_S)$).
\end{proposition}
\begin{proof}
As above, the argument is similar in both cases so we only treat the $\R$-constructible case. We may choose as above a $\mu$-stratification $X=\bigsqcup_{\alpha} X_{\alpha}$ such that $G|_{X_{\alpha}\times T}$ is isomorphic to $p_{X_{\alpha}}^{-1}G_{\alpha}$ for some $G_{\alpha}\in\rD^\rb_{\coh}(\sho_T)$. The assumption entails that $\pi$ is proper on the support of each $G_{\alpha}$. Since $(R\pi_*G)|_{X_{\alpha}\times S}=R\pi_{\alpha *}(G|_{X_{\alpha}\times T})$ where $\pi_{\alpha}:=\pi|_{X_{\alpha}\times T}$, the result follows by Grauert's direct image theorem (\cf \cite{HG}).
\end{proof}

When $X$ is a complex analytic manifold, we  refer to \cite{MFCS1} and \cite{MFCS2} for details on relative holonomy for relative $\DXS$-modules. We recall that a holonomic $\DXS$-module $\shm$ is regular if the holomorphic restrictions of $\shm$ to the fibers of $p_X$ are regular holonomic $\shd_X$-modules. We note by $\rD^\rb_{\rhol}(\DXS)$ the bounded triangulated category of complexes with regular holonomic cohomology. We note by $d_X$ (resp. by $d_S$) the complex dimension of $X$ (resp. of $S$).

\begin{itemize}
\item
{
We recall that $\bD_X\pSol(\cbbullet)$ is functorially isomorphic to $
\pDR_X(\cbbullet):=\Rhom_{\DXS}(\sho_{X\times S}, \cbbullet)[d_X]$.}
\item{Following (3.1) in \cite{MFP2}, we note by $\rho_S: X\times S\to X_{sa}\times S$  the natural morphism of sites where $X_{sa}$ denotes the subanalytic site underlying $X$. We note by $\sho^{t,S}_{X\times S}$ the relative subanalytic sheaf on $X_{sa}\times S$ associated to the subanalytic sheaf $\sho^{t}_{X\times S}$ on $(X\times S)_{sa}$ introduced in \cite{KS5}.
Then, for $F\in\rD^\rb_{\rc}(\pOS)$, one defines
$$\RH^S_X(F):=\rho_S^{-1}\Rhom_{\rho_{S*}\pOS}(\rho_{S*}F, \sho^{t,S}_{X\times S})[-d_X]$$ which is an object of $\rD^\rb(\DXS)$.
}
\item{As proved in \cite[Th.~1]{FMFS1} \cite[Th.~2]{FMFS2} to which we refer for details, the category $\rD^\rb_{\cc}(\pOS)$ is equivalent to the category $\rD^\rb_{\rhol}(\DXS)$ via the functor $\RH^S_X$ and its quasi-inverse $\pSol_X:=\Rhom_{\DXS}(\cbbullet, \sho_{X\times S})[d_X]$.}
\end{itemize}

\begin{proposition}\label{Linvim}
Let $\pi: T\to S$ be a morphism of complex manifolds and let $X$ be given. Then there exists a natural transformation of functors  $\Psi_{\cbbullet}$ from $\rD^\rb_{\rc}(\pOS)$ to $\rD^\rb(\shd_{X\times T/T})$ $$\Psi_F: L\pi^*\RH^S_X(F)\to \RH^{T}_X(L\pi^*F).$$ If, moreover, $F\in\rD^\rb_{\cc}(\pOS)$ then $\Psi_F$ is an isomorphism in $\rD^\rb_{\rhol}(\shd_{X\times T/T})$. 
\end{proposition}
\begin{proof} 

We refer to \cite{KS5} and \cite{KS6} for the background on sites and Grothendieck 
topologies. We start by constructing a natural morphism in $\rD^\rb(\shd_{X\times T/T})$
$$\Psi_F:L\pi^*\RH^S_X(F)\to \RH^{T}_X(L\pi^*F)$$
for $F\in\rD^\rb_{\rc}(p_X^{-1}\sho_S)$.

 Following the definition of the relative Riemann-Hilbert functor, we aim to construct a natural
morphism:
\begin{multline}\label{E3}
\sho_{X\times T}\otimes^L_{\pi^{-1}\sho_{X\times S}}\pi^{-1}\RH^S_X(F)=\\
\sho_{X\times T}\otimes^L_{\pi^{-1}\sho_{X\times S}}\pi^{-1}{{\rho}_{S}^{ -1}}\Rhom_{\rho_{S *}p^{-1}\sho_S}(\rho_{{S} *}F,\sho^{t,S}_{X\times S})[-d_X]\\
\to\rho^{-1}_{T}\Rhom_{\rho_{T *}p^{-1}\sho_{T}}(\rho_{ T *}L\pi^*F, \sho^{t,T}_{X\times T})[-d_X]=\RH^T_X(L\pi^*F).
\end{multline}
The isomorphism of functors on sites $\pi^{-1}\rho^{-1}_S\simeq \rho^{-1}_{T}\pi^{-1}$ yields a natural isomorphism in $\rD^\rb(\pi^{-1}\DXS)$
\begin{multline}\label{E4}
\pi^{-1}{\rho^{-1}_{S}}\Rhom_{\rho_{S *}p^{-1}\sho_S}(\rho_{S *}F,\sho^{t,S}_{\XS})\\
\simeq {\rho^{-1}_{T}}\pi^{-1}\Rhom_{\rho_{S *}p^{-1}\sho_S}(\rho_{S *}F,\sho^{t,S}_{\XS}).
\end{multline}
Recall that, for a morphism $f: Z'\to Z$ of manifolds and $\sha$ a sheaf of rings on $Z$, we have a natural morphism of bifunctors on $\rD^\rb(\sha)$ (\cf \cite[(2.6.27]{KS1}):
\begin{equation}\label{eq:fA}
f^{-1}\Rhom_{\sha}(\cbbullet, \cbbullet)\to \Rhom_{f^{-1}\sha}(f^{-1}(\cbbullet), f^{-1}(\cbbullet)).
\end{equation}
Since we are working with sheaves on Grothendieck topologies, we have the analogous of \eqref{eq:fA}, that is, we have a natural morphism
in $\rD^\rb(\rho_{T !}\pi^{-1}\DXS)$
\begin{equation}\label{E70}
\pi^{-1}\!\!\Rhom_{\rho_{ S *}p^{-1}\sho_S}(\rho_{S *}F,\sho^{t,S}_{\XS})
\to \Rhom_{\pi^{-1}\rho_{S *}p^{-1}\sho_S}(\pi^{-1}\rho_{S *}F,\pi^{-1}\sho^{t,S}_{\XS}),
\end{equation}
hence a natural morphism
in $\rD^\rb(\pi^{-1}\DXS)$
\begin{multline}\label{E5}
{\rho^{-1}_{T}}\pi^{-1}\Rhom_{\rho_{S *}p^{-1}\sho_S}(\rho_{S *}F,\sho^{t,S}_{\XS})\\
\to {\rho^{-1}_{T}}\Rhom_{\pi^{-1}\rho_{S *}p^{-1}\sho_S}(\pi^{-1}\rho_{S *}F,\pi^{-1}\sho^{t,S}_{\XS}).
\end{multline}
Thus we get a functorial chain of morphisms in $\rD^\rb(\pi^{-1}\DXS)$
\begin{align*}
\rho^{-1}_{T}\pi^{-1}&\Rhom_{\rho_{S *}p^{-1}\sho_S}(\rho_{S *}F,\sho^{t,S}_{\XS})\\
&\to \rho^{-1}_{T}\Rhom_{\pi^{-1}\rho_{ S *}p^{-1}\sho_S}(\pi^{-1}\rho_{S *}F,\sho^{t,T}_{X\times T})\\
%\hspace*{2.3cm}
&\simeq \rho^{-1}_{T}\Rhom_{\rho_{T *}p^{-1}\sho_{T}}(\rho_{T*}p^{-1}\sho_{T}\otimes^L _{\pi^{-1}\rho_{S*}p^{-1}\sho_S}\pi^{-1}\rho_{S*}F,\sho^{t,T}_{X\times T})\\
&\simeq \rho^{-1}_{T}\Rhom_{\rho_{T*}p^{-1}\sho_{T}}(\rho_{T*}L\pi^*F,\sho^{t,T}_{X\times T}),
\end{align*}
where the last isomorphism follows by \cite[Lem. 5.7]{FMFS2}.
We also
remark that the last term shifted by $[-d_X]$
(which we will name $\shl$ for simplicity) 
is
the right term of the desired morphism \eqref{E3} and is already an object of $\rD^{\rb}(\shd_{X\times T/T})$.
Hence, by applying to \eqref{E4} the derived functor $\sho_{X\times T}\otimes^L_{\pi^{-1}\sho_{\XS}}(\cbbullet)[-d_X]$, we obtain a chain of natural morphisms
in $\rD^\rb(\shd_{X\times T/T})$.
\begin{multline}\label{E6}
\sho_{X\times T}\otimes^L_{\pi^{-1}\sho_{\XS}}\pi^{-1}\rho^{-1}_{S}\Rhom_{\rho_{S *}p^{-1}\sho_S}(\rho_{S *}F,\sho^{t,S}_{\XS})[-d_X]\\
\to \sho_{X\times T}\otimes^L_{\pi^{-1}\sho_{\XS}}\shl\to\shl
\end{multline}
The composition of the morphisms in $(\ref{E6})$ gives the desired morphism $\Psi_F$.

If now we assume that $F\in\rD^\rb_{\cc}(\pOS)$, according to \cite[Prop.~3.10]{FMFS2} and Proposition~\ref{PN2}, $L\pi^*\RH_X^S(F)$ and $\RH_X^{T}(L\pi^*F)$ are complexes with regular holonomic cohomology. We are thus in conditions of applying the variant of Nakayama's lemma (\cf \cite[Cor.~1.10]{MFCS2}). That is, to end the proof it is sufficient to prove that,  after composing  $Li^*_t$ with $\Psi_F$, we get un isomorphism for any $t\in T$ which follows by \cite[Prop.~5.9]{FMFS2}. 
\end{proof}

Let us introduce the full triangulated subcategory $\rD^\rb_{\rhol,\rm{good}}(\shd_{X\times S/S})$ 
of $\rD^\rb_{\rhol}(\shd_{X\times S/S})$ whose objects are complexes with $\sho_{X\times S}$-good cohomology\footnote{Similarly to the absolute case, when $d_X=d_S=1$ any relative regular holonomic 
$\DXS$-module is good. However we do not know if this is true in general.}.

For a morphism $f:Y\to X$ or a projective morphism $\pi:T\to S$, 
following \cite[Def.~ 2.10]{FMFS2},
let us introduce the full triangulated subcategory $\rD^\rb_{\rhol,f\good}(\shd_{Y\times S/S})$ of $\rD^\rb_{\rhol}(\shd_{Y\times S/S})$ (respectively the subcategory $\rD^\rb_{\rhol,\pi\good}(\shd_{X\times T/T})$ of the category $\rD^\rb_{\rhol}(\shd_{X\times T/T})$) whose objects are complexes with $f$-good cohomology (respectively with $\pi$-good cohomology).
\begin{proposition}\label{Ldirectim}
Assume that $\pi:T\to S$ is a projective morphism and let $\Omega_{T/S}$ denote the sheaf of relative forms of maximal degree. Then, for any $G\in \rD^\rb_{\cc}(p^{-1}\sho_T)$ such that $\RH^T_X(G)$ belongs to $\rD^\rb_{\rhol,\pi-\rm{good}}(\shd_{X\times T/T})$, we have a natural isomorphism  $$R\pi_*\RH^T_X(G)[d_T]\to \RH^S_X(R\pi_*(G\otimes_{p^{-1}\sho_{T}}
p^{-1}\Omega_{T/S}))[d_S]$$
which gives 
$R\pi_*\RH^T_X(G)\simeq  
\RH^S_X(R\pi_*(G\otimes_{p^{-1}\sho_{T}}p^{-1}(\pi^!\sho_S)))
$ (since $\pi^!\sho_S=\Omega_{T/S}[d_T-d_S]$).
\end{proposition}
\begin{proof}
According to \cite[Th.~3.15]{FMFS2}, by the relative Riemann-Hilbert correspondence it is sufficient to prove that we have a natural isomorphism after applying the De Rham functor to both sides.
For any $\shm\in\rD^\rb_{\rhol}(\shd_{X\times T/T})$
$R\pi_*\DR(\shm)\overset{\simeq}{\to} \DR(R\pi_*\shm)$, 
which follows by adjunction:
\begin{align*}R\pi_*(\Rhom_{\shd_{X\times T/T}}(\sho_{X\times T}, \shm))&\simeq R\pi_*(\Rhom_{\pi^{-1}\shd_{X\times S/S}}(\pi^{-1}\sho_{X\times S}, \shm))\\
&\simeq \Rhom_{\DXS}(\sho_{X\times S}, R\pi_*\shm).\end{align*}
Since $\DR\circ \RH_X^S(\cbbullet)\simeq \bD'(\cbbullet)$ the desired morphism 
is given by
\begin{align*}
R\pi_*\bD'(G)[d_T] =&R\pi_*\Rhom_{p^{-1}\sho_{T}}(G,p^{-1}\sho_{T})[d_T]\\
\to &
R\pi_*\Rhom_{p^{-1}\sho_{T}}(p^{-1}\Omega_{T/S}\otimes_{p^{-1}\sho_{T}}G,
p^{-1}\Omega_{T/S})[d_T]\\
\to& \Rhom_{p^{-1}\sho_{S}}(R\pi_*(p^{-1}\Omega_{T/S}\otimes_{p^{-1}\sho_{T}}G),
R\pi_*(p^{-1}\Omega_{T/S}))[d_T]\\
\to &\Rhom_{p^{-1}\sho_{S}}(R\pi_*(p^{-1}\Omega_{T/S}\otimes_{p^{-1}\sho_{T}}G),
p^{-1}\sho_S)[d_S]= \\
=&
\bD'(R\pi_*(p^{-1}\Omega_{T/S}\otimes_{p^{-1}\sho_{T}}G))[d_S]\\
\end{align*} 
where the last morphism is obtained by composing with the trace map
$R\pi_*\pi^!\sho_S\to \sho_S$ .
Now let us prove it is an isomorphism.
We can choose a $\mu$-stratification $X=\bigsqcup_{\alpha} X_{\alpha}$ such that $G|_{X_{\alpha}\times T}$ is isomorphic to $p_{X_{\alpha}}^{-1}G_{\alpha}$ for some $G_{\alpha}\in\rD^\rb_{\coh}(\sho_T)$.
The previous morphism  is an isomorphism since 
%left term of  restricted to $X_\alpha\timeq S$ is
\begin{align*}
R\pi_*\bD'(G)[d_T]|_{X_{\alpha}\times S}&=
p_{X_\alpha}^{-1}R\pi_*\bD'(G_\alpha)[d_T]\simeq p_{X_\alpha}^{-1}\bD'(R\pi_*(\Omega_{T/S}\otimes_{\sho_{T}}G_\alpha))[d_S]\\
&=
\bD'(R\pi_*(p^{-1}\Omega_{T/S}\otimes_{p^{-1}\sho_{T}}G))[d_S]|_{X_{\alpha}\times S}\\ \end{align*}
by Serre's duality for complex manifolds.
% and Lemma \ref{L1} below.
\end{proof}

The following lemma can be easily deduced using the proof of \cite[Lem.~3.22]{MFCS2}.

\begin{lemma}\label{LRFM}
Let G be a coherent $\sho_S$-module and let $F\in\rD^\rb_{\rc}(\pOS)$. Then we have an isomorphism in $\rD^\rb_{\rhol}(\DXS)$ 

$\Rhom_{\pOS}(p_X^{-1}G, \RH_X^S(F))\simeq \RH^S_X(p_X^{-1}G\otimes^L_{\pOS}F)$ functorial in $F$. In particular
\begin{itemize}
\item{If  $\RH_X^S(F)$ is a good-$\DXS$-module, then $\RH^S_X(p_X^{-1}G\otimes^L_{\pOS}F)$ is also good.}
\item{$\pSol(\Rhom_{\pOS}(p_X^{-1}G, \RH_X^S(F)))\simeq p_X^{-1}G\otimes^L_{\pOS}F.$}
\end{itemize}
\end{lemma}

\subsection{Relative constructibility when $S$ is projective}\label{pre2}We recall the following definition in \cite[2.1]{BonVdBergh}:

\begin{definition}
Given a set of objets ${\mathcal E} = (E_i)_{i\in I}$ of a triangulated category $\shd$, we denote by $\langle \mathcal E\rangle$ the smallest epaisse (thick) triangulated full subcategory of $\shd$ which contains $\she$ (the so called epaisse envelope of $\she$ which is thus closed under shifts, isomorphisms, direct summands and distinguished triangles). We say that $\she$ classically generates a triangulated category $\mathcal D$ if $\mathcal D=\langle \mathcal E\rangle$. 
\end{definition}

\begin{remark}\label{RemEp} 
Let us suppose that $\she$ classically generates a triangulated category $\mathcal D$.
Let $\mathcal T$ be a property on $\mathcal D$ such that
\begin{enumerate}
\item $X_i\in \shd$ satisfy $\mathcal T$ for  $i=1,2$ if and only if  $X_1\oplus X_2$ satisfies $\mathcal T$;
\item if $X\to Y\to Z\to X[1]$ is a distinguished triangle such that $X$ and $Y$ satisfy $\mathcal T$ then $Z$ satisfies $\mathcal T$;
\item for any $E\in \she$ and for any $n\in \Bbb Z$ the object $E[n]$ satisfies $\mathcal T$.
\end{enumerate}
Then any object of $\shd$ satisfies $\sht$.
In fact the previous conditions entail that the full subcategory $\shd'$ of $\shd$ whose objects satisfy 
$\sht$ is a strictly full (that is, full and closed by isomorphisms) thick triangulated subcategory of $\shd$ and hence $\shd=\langle \she \rangle\subseteq \shd'$. 
We see that $0$ satisfies property $\sht$ since given $X\in\she$  the distinguished triangle $X\overset{\id}{\to}X\to 0\overset{+1}{\to}$ implies $0\in\shd'$. 
This implies that $\shd'$ is closed under isomorphisms.
\end{remark}

Let us introduce the following condition:
\begin{assumption}\label{Assumption}
We assume that $S$ is  a projective complex algebraic manifold, thus $S$ is the analytification of a smooth proper complex variety $S^{\mathrm{alg}}$.
\end{assumption}

In the remaining of this section we suppose that $S$ satisfies Assumption~\ref{Assumption}. 
 The  result below is known by the specialists but, since we have not found a reference, we give a detailed proof.

\begin{proposition}\label{Remalgan}
The functor $E$ in \eqref{E200}
is an equivalence of categories.
\end{proposition}
\begin{proof}
Let us consider the following commutative square of functors:
\begin{equation}\label{E300}
\xymatrix{
\rD^\rb(\coh(\sho_{S^{\mathrm{alg}}}))\ar[d]\ar[r] & \rD^\rb_{\coh}(\sho_{S^{\mathrm{alg}}})\ar[d]\\
\rD^\rb(\coh(\sho_S))\ar[r]^E & \rD^\rb_{\coh}(\sho_S)\\
}
\end{equation}
We have to prove that the bottom arrow is an equivalence of categories.
By induction on the cohomological length of complexes (using 
 truncated triangles) we get the essential  surjectivity. 
 Let us prove that $E$ is fully faithful.
It is well known (for instance, see \cite{BN}) that the top horizontal arrow
is an equivalence.
By Serre's GAGA Theorem the left vertical arrow (given by the analytification functor )
is an equivalence. 
Hence $E$ is fully faithful if and only if the right vertical arrow in \eqref{E300} is fully faithful.
Let us recall that any coherent $\sho$-module on a smooth projective complex algebraic variety  is the epimorphic image of a locally free $\sho$-module of finite rank, hence, due to Assumption~\ref{Assumption}, we have $\rD^\rb_\coh (\sho_{S^{\mathrm{alg}}})\simeq \rD^\rb(\coh (\sho_{S^{\mathrm{alg}}}))$.
 
Given $\shl_1, \shl_2$ locally free $\sho_{S^{\mathrm{alg}}}$-modules of finite rank, we have 
the following isomorphisms  where in the fourth and sixth we use Serre's GAGA Theorem:
\begin{align*}
\Hom_{\rD^\rb_{\coh}(\sho_{S^{\mathrm{alg}}})}&(\shl_1, \shl_2[n])=
\Hom_{\rD^\rb(\sho_{S^{\mathrm{alg}}})}(\shl_1, \shl_2[n])\\
=&
\Hom_{\rD^\rb(\sho_{S^{\mathrm{alg}}})}(\sho_{S^{\mathrm{alg}}}, \bD\shl_1\otimes_{\sho_{S^{\mathrm{alg}}}}\shl_2[n])\\
=&{\mathrm H}^n(S^{\mathrm{alg}}, \bD\shl_1\otimes_{\sho_{S^{\mathrm{alg}}}}\shl_2)=
{\mathrm H}^n(S, (\bD\shl_1\otimes_{\sho_{S^{\mathrm{alg}}}}\shl_2)^{\mathrm{an}})\\
=&
{\mathrm H}^n(S, (\bD\shl_1)^{\mathrm{an}}\otimes_{\sho_{S}}\shl_2^{\mathrm{an}})=
{\mathrm H}^n(S, \bD\shl_1^{\mathrm{an}}\otimes_{\sho_{S}}\shl^{\mathrm{an}}_2)\\
=&\Hom_{\rD^\rb(\sho_{S})}(\shl^{\mathrm{an}}_1, \shl^{\mathrm{an}}_2[n])
=\Hom_{\rD^\rb_{\coh}(\sho_{S})}(\shl^{\mathrm{an}}_1, \shl^{\mathrm{an}}_2[n]).\\
\end{align*}

Let us fix  a locally free $\sho_{S^{\mathrm{alg}}}$-module $\shl_1$ of finite rank.
We will say that an object $F\in \rD^\rb_{\coh}(\sho_{S^{\mathrm{alg}}})$ satisfies the property
 $\mathcal T$ if the canonical morphism
\begin{equation}\label{EqFF}
\Hom_{\rD^\rb_{\coh}(\sho_{S^{\mathrm{alg}}})}(\shl _1,F )\to
 \Hom_{\rD^\rb_{\coh}(\sho_{S})}(\shl_1^{\mathrm{an}},F^{\mathrm{an}})
 \end{equation} is an isomorphism.
For this property we just proved condition $(3)$ of  Remark~\ref{RemEp} and condition $(1)$ is clear. 
If $F_1\to F_2\to F_3\to F_1[1]$ is a distinguished triangle 
by \cite[1.1.1]{BBD} we have the following commutative diagram whose rows are exact
{\footnotesize{\[
\xymatrix@-20pt{
\dots\ar[r]& \Hom_{\rD^\rb_{\coh}(\sho_{S^{\mathrm{alg}}})}(\shl_1,F_1)
\ar[r] \ar[d]&\Hom_{\rD^\rb_{\coh}(\sho_{S^{\mathrm{alg}}})}(\shl_1,F_2 )\ar[r]\ar[d]  &
\Hom_{\rD^\rb_{\coh}(\sho_{S^{\mathrm{alg}}})}(\shl_1,F_3 )\ar[r] \ar[d] 
& \dots \\  
\dots\ar[r] & \Hom_{\rD^\rb_{\coh}(\sho_{S})}(\shl^{\mathrm{an}}_1,F^{\mathrm{an}}_1 )
\ar[r] &\Hom_{\rD^\rb_{\coh}(\sho_S)}(\shl^{\mathrm{an}}_1,F^{\mathrm{an}}_2 )\ar[r]
&
\Hom_{\rD^\rb_{\coh}(\sho_S)}(\shl^{\mathrm{an}}_1,F^{\mathrm{an}}_3 )\ar[r]  
&\dots \\  }\]}}so condition $(2)$ follows by the five Lemma.
This proves that \eqref{EqFF} is an isomorphism for any $F\in \rD^\rb_{\coh}(\sho_{S^{\mathrm{alg}}})$
and for any $\shl_1\in\she$.
By the same argument in the first variable, we get that
$\Hom_{\rD^\rb_{\coh}(\sho_{S^{\mathrm{alg}}})}(G,F )\to
 \Hom_{\rD^\rb_{\coh}(\sho_{S})}(G^{\mathrm{an}},F^{\mathrm{an}})$ is an isomorphism
 for any $F,G\in \rD^\rb_{\coh}(\sho_{S^{\mathrm{alg}}})$.
 This concludes the proof that the right arrow in \eqref{E300} is fully faithful.
 \end{proof}

\begin{remark}\label{Rmkfree}
Given a complex $F\in \rD^\rb(\coh(\sho_S))$
there exists a complex of locally free $\sho_S$-modules of finite rank $L(F)$ and a 
a quasi-isomorphism  $L(F)\to F$; thus, by the previous result, for any $G$ in $\rD^\rb_{\coh}(\sho_S)$, there exists a complex of locally free $\sho_S$-modules of finite rank $L(G)$ and 
an isomorphism  $L(G)\to G$ in $\rD^\rb_{\coh}(\sho_S)$.
Moreover the set  of locally free sheaves of finite rank classically generates 
$\rD^\rb_{\coh}(\sho_S)$.
\end{remark}

\begin{lemma}\label{L1}
Let ${\mathcal E}$ be the set of objects  $H$ in $\rD^\rb_{\rc}(\pOS)$ which are (locally finite) direct sums of terms of the form $\C_{U}\boxtimes L$ where $U$ is relatively compact contractible subanalytic open in $X$
and $L$ is a 
locally free $\sho_S$-module of finite rank. Then $\she$ classically generates $\rD^\rb_{\rc}(\pOS)$. 
\end{lemma}
\begin{proof}
We aim to prove that  $\rD^\rb_{\rc}(\pOS)=\langle \mathcal E\rangle$. Because the right side is contained in the left side, we are led to prove that, given $F\in \rD^\rb_{\rc}(\pOS)$, then $F$ belongs to $\langle \mathcal E\rangle$.
We will use the notations of \cite[\S8.1]{KS1}. Let us consider a subanalytic stratification $(X_{\alpha})_{\alpha\in A}$ of $X$ such that $F|_{X_{\alpha}\times S}$ is locally isomorphic to $p^{-1}_{X_{\alpha}}G$, where $G\in\rD^\rb_{\coh}(\sho_S)$. We then choose a subanalytic triangulation of $X$ adapted to $(X_{\alpha})_{\alpha\in A}$ provided by theorem (\cite[Th.\,8.2.5]{KS1}). Let $\textbf{T}=(T,\Delta)$ be the associated simplicial complex  on $X$. We identify each simplex to its homeomorphic image in $X$, so that each $X_{\alpha}$ is a union of simplices and each simplex is contained in some $X_{\alpha}$. 

For any simplex $\sigma$, its open star $U(\sigma)$ is contractible. The argument already
used in \cite[Proof of Prop.\,3.3(i))]{MFCS2} shows that the entries of $F$ are $p_{U(\sigma)*}$-acyclic and that, for any $x\in|\sigma|$, for any open subset $V$ of $S$ 
\begin{equation}\label{eq:VFx}
\Gamma(U(\sigma)\times V; F)\isom \Gamma(V; F|_{\{x\}\times S}).
\end{equation} is an isomorphism.

We conclude that the natural morphism of complexes
\begin{equation}\label{E0}
G_{\sigma}(F):=p_{U(\sigma) \ast} F|_{U(\sigma)\times S}\to F|_{\{x\}\times S},
\end{equation}
is a quasi-isomorphism for any $x\in |\sigma|$. 
By~adjunction there exists a canonical morphism
\begin{equation}\label{E01}
p_{U(\sigma)}^{-1}G_\sigma(F)\to F|_{U(\sigma)\times S},
\end{equation}
and by extension a canonical morphism of complexes
\begin{equation}\label{1}
u_{G_\sigma(F)}: \C_{U(\sigma)\times S}\otimes p_X^{-1}G_\sigma(F)\to F.
\end{equation}
Thus, for any $\sigma$, for any $x\in |\sigma|$, and any $s\in S$, the germ
\begin{equation}\label{E002}
u_{G_\sigma(F),\,(x, s)}: (\C_{U(\sigma)}\boxtimes G_\sigma( F))_{(x,s)}\to F_{(x,s)}
\end{equation}
is a quasi-isomorphism.

Because $F|_{\{x\}\times S}\in\rD^\rb_{\coh}(\sho_S)$,   
 by Remark~\ref{Rmkfree}
there exists  an isomorphism $L_\sigma(F)\to G_\sigma(F)$ in $\rD^\rb_{\coh}(\sho_S)$
with $L_\sigma(F)$
a bounded complex  
 of locally free $\sho_S$-modules of finite rank.
 
 By \eqref{E01} we obtain a morphism in the derived category 
 \begin{equation}\label{2}
p_{U(\sigma)}^{-1}L_\sigma(F)\to F|_{U(\sigma)\times S}
\end{equation}
and thus a natural morphism in $\rD^\rb_{\rc}(\pOS)$
\begin{equation}\label{E02}
u_{L_\sigma(F)}: \C_{U(\sigma)\times S}\otimes p_X^{-1}L_\sigma(F)\to F
\end{equation} 
such that
\begin{equation}\label{E02f}
u_{L_\sigma(F),\,(x, s)}: (\C_{U(\sigma)}\boxtimes L_\sigma( F))_{(x,s)}\to F_{(x,s)}
\end{equation} 
is an isomorphism, for any $\sigma$, for any $x\in |\sigma|$, and any $s\in S$.

We shall now adapt the proof of Proposition A.2 (i) of \cite{D-G-S11}. 
We shall prove by induction on $i$ that there exists a morphism $u_i: G_i\to F$ in $ \rD^\rb_{\rc}(\pOS)$  such that:

\begin{enumerate}
\item\label{ast1}
each $G_i$ belongs to $\langle\mathcal E\rangle$,
\item\label{ast2} $u_{i|\;|\bf{T}_i|\times S}: G_{i|\;|\bf{T}_i|\times S}\to F_{|\;|\bf{T}_i|\times S}$ is an isomorphism.
\end{enumerate}
The desired result follows with $i=d_X$.

We start by noticing that, if two terms in a distinguished triangle satisfy \eqref{ast1}, then the the same is true for the third one.
For $i=0$ we set $
u_0=\oplus_{\sigma\in |\bf{T}_0|}  u_{L_\sigma(F)}$ which satisfies the desired condition in view of \eqref{E02f}.
Let us suppose that we have constructed a morphism  $u_i: G_i\to F$ in $ \rD^\rb_{\rc}(\pOS)$ 
  satisfying \eqref{ast1} and \eqref{ast2}.
Let us complete   $u_i$ into a distinguished triangle
$H_i \stackrel{v_i}\to G_i\stackrel{u_i}\to F \stackrel{+}\to $.
By construction $H_i|_{|\bf{T}_i|\times S}$ is isomorphic to $0$ which implies that
$\oplus_{\sigma \in\Delta_{i+1}\setminus \Delta_i}(H_i)_{|\sigma |\times S}\to H_i|_{|\bf{T}_{i+1}|\times S}$ is an isomorphism.

Because $H_i$ is a complex in $\rD^\rb_{\rc}(\pOS)$, by \eqref{E02} we can construct a morphism 
$$u'_{i+1}:=\oplus_{\sigma\in \Delta_{i+1}\setminus \Delta_i}  u_{L_\sigma(H_i)}: G'_{i+1}:=\oplus_{\sigma\in \Delta_{i+1}\setminus \Delta_i}\C_{U(\sigma)\times S}\otimes p_X^{-1}L_\sigma(H_i)\to H_i$$
in $ \rD^\rb_{\rc}(\pOS)$. In view of the remark above, the restriction of the latter to $|\bf{T}_{i+1}|$ is an  isomorphism since by \eqref{E02f}
$(u'_{i+1})_{(x,s)}$ is an  isomorphism for any $x\in |\bf{T}_{i+1}|\setminus |\bf{T}_i|$ and for any $s\in S$.

Let us choose  $H_{i+1}$ and $G_{i+1}$ in $\rD^\rb_{\rc}(\pOS)$ such that $u'_{i+1} $
and $v_i\circ u'_{i+1}$ insert into distinguished triangles
$$G'_{i+1}\overset{u'_{i+1}}{\to}H_i\to H_{i+1}\overset{+1}{\to}$$
$$G'_{i+1}\overset {v_i\circ u'_{i+1}}{\to} G_i\to G_{i+1}\overset{+1}{\to}$$
By the definition of $G'_{i+1}$ and the condition on $G_i$ given by the induction hypothesis, we see that  $G_{i+1}$ satisfies condition \eqref{ast1}. 

As in 
 \cite[Prop.~A.2 (i)]{D-G-S11} the octahedral axiom applied to the preceding triangles gives a morphism
 $u_{i+1}:G_{i+1}\to F$ and a distinguished triangle
 $$H_{i+1}\to G_{i+1}\overset{u_{i+1}}{\to} F\overset{+1}{\to}$$ 
 in $\rD^\rb_{\rc}(\pOS)$.
 
By construction $H_{i+1}|_{|\bf{T}_{i+1}|}$ is isomorphic to $0$. 
Thus $u_{i+1}$ satisfies \eqref{ast1} and  \eqref{ast2}.
\end{proof}

 \remark
The following Theorem  leads to a wide simplification in the proof of the relative Riemann-Hilbert correspondence (\cite[Th.~2]{FMFS2}) if one assumes Assumption~\ref{Assumption}.
 Moreover it seems expectable that such a result would be essential for a generalization of the relative Riemann-Hilbert correspondence to the case of a projective smooth morphism $X\to S$.

\begin{theorem}\label{Cor1}
The natural  functor
\begin{equation}\label{eqcoh}
\shi:\rD^\rb(\Mod_{\rc}(\pOS))\to \rD^\rb_{\rc}(\pOS)\end{equation}
 is an equivalence of categories.
\end{theorem}
\begin{proof}

In order to prove that  
 $\shi$ is fully faithful 
it is sufficient to prove that,
given $F, G \in \rD^\rb(\Mod_{\rc}(\pOS))$, the morphism
\begin{equation}\label{FF}
\shi_{F,G}:\Hom_{\rD^\rb(\Mod_{\rc}(\pOS))}(F,G)\to \Hom_{\rD^\rb(\Mod_{\wrc}(\pOS))}(F,G)
\end{equation}
is bijective. 

Let us prove that $\shi_{F,G}$ is injective. 
We start by noticing that if a complex $F$ belongs to $\rD^\rb(\Mod_{\rc}(\pOS))$ then, for any $s\in S$, $Li^*_sF$ is a complex in $\rD^\rb(\Mod_{\rc}(X))$.
Let $\phi\in \Hom_{\rD^\rb(\Mod_{\rc}(\pOS))}(F,G)$ be represented by the diagram below in $\rK^\rb(\Mod_{\rc}(\pOS))$
$$
\xymatrix@R=6pt{
F & & G.\\
& F'\ar[ul]_{\scriptstyle{{\mathrm{qis}}} }\ar[ru]^{\phi} & \\
}
$$
If $\shi_{F,G}(\phi)=0$ as a morphism in $\Hom_{\rD^\rb_{\rc}(\pOS)}(F,G)$,  we can insert the preceding diagram in the following commutative one, where $E\in\rK^\rb(\Mod_{\wrc}(\pOS))$, 
$$
\xymatrix@R=8pt{
F & & G.\\
& F'\ar[ul]_{\scriptstyle{{\mathrm{qis}}} }\ar[ru]^{\phi} & \\
& E\ar@/^0.4pc/[uul]^(0.3){\scriptstyle{{\mathrm{qis}}} }\ar[u]^(0.6){\scriptstyle{{\mathrm{qis}}} } \ar@/_0.4pc/[ruu]_(0.3){0} & \\
& \\
}
$$
Thus we have $\phi\circ \mathrm{qis}=0$. For any $s\in S$, studying the commutative diagram which is obtained by applying the functor $Li^*_s$ to the above one, we get $Li_s^ *\phi\circ Li_s^*\mathrm{qis}=0$. Since $\Hom_{\rD^\rb(\Mod_{\rc}(X))}(Li^*_sF, Li^*_sG)\to  \Hom_{\rD^\rb_{\rc}(X)}(Li^*_sF, Li^*_sG)$ is an isomorphism, we conclude that $Li^*_s\phi=0$ thus, by the relative variant of Nakayama's Lemma, we get $\phi=0$.

Let us  now prove the surjectivity of \eqref{FF}.
Let us consider the following cartesian  diagram 
\[ 
\xymatrix{
X\times S \ar[r]^{p}\ar[d]_{q} \ar[rd]^{f}& S\ar[d]^{a} \\
X\ar[r]^{b}& \{pt\}}
\]
Let $\mathcal E$ be as in Lemma \ref{L1}. 

Let us first assume that $\shl=\C_{U}\boxtimes L$ of $\she$. We have: 
\begin{align*}
\R \Gamma (X\times S; \shl)=&
R f_*(\shl)=
\R b_*\R q_*(q^{-1}\C_{U}\otimes p^{-1}L)\\
=&\R b_*(\C_{U}\otimes \R q_* p^{-1}L)\\
=&\R b_*(\C_{U}\otimes b^{-1}\R a_* L)\\
=&\R b_*\C_{U}\otimes\R a_*L\\
=&\R \Gamma (X; \C_{U})\otimes \R \Gamma (S; L)
\\
\end{align*}
where the third and the 
fifth equation are given by the projection formula.

In view of Remark \ref{RemEp} and similarly to the proof of Proposition~\ref{Remalgan} it is enough to prove that
 \eqref{FF} is an isomorphism (we already proved it is injective) when $\shf, \shg$ belong to $\mathcal E$ up to a shift. Next, by the definition of $\mathcal E$, it is sufficient to assume that $\shf=\C_{U}\boxtimes F$ 
and $\shg=\C_{V}\boxtimes G$.

We have:
\begin{align*}
\Hom_{\rD^\rb_{\rc}(\pOS)}&(\shf, \shg[k])\\=&\Hom_{\rD^\rb(\pOS)}(\shf,\shg[k])\\=&\Hom_{\rD^\rb(\puOS)}(p_U^{-1}F,\shg|_{U\times S})[k]
\\=&\Hom_{\rD^\rb(\puOS)}(p_U^{-1}\sho_S,\C_{V}
\boxtimes (F^\vee\otimes_{\sho_S}G))[k]\\
=&{\mathrm H}^k(U\times S,\C_{V}\boxtimes (F^\vee\otimes_{\sho_S}G))
\\=&\bigoplus\limits_{p+q=k}
{\mathrm H}^p(U,\C_{V})\otimes_{\Bbb C}{\mathrm H}^q(S, F^\vee\otimes_{\sho_S}G)\\
\end{align*}
where the last equation is obtained by taking $H^k$ of 
$\R \Gamma (X\times S; \C_{U}\boxtimes L)=\R \Gamma (X; \C_{U})\otimes \R \Gamma (S; L)$.
 Let us consider the following commutative diagram, where, by our assumption, the bottom left vertical arrow is an isomorphism.
{\tiny{\[
\xymatrixcolsep{0.3in}\xymatrixrowsep{3.5em}
\xymatrix@-15pt{
\bigoplus\limits_{\substack{\scriptscriptstyle{p+q=k}}}{\mathrm H}^p(U,\C_{V})\otimes_{\Bbb C}{\mathrm H}^q(S, F^\vee\otimes_{\sho_S}G)
\ar[r]^(0.55)*[@]{\cong}
\ar[d]^*[@]{\cong}&  \Hom_{\rD^\rb(\pOS)}(\shf, \shg[k])\\
\bigoplus\limits_{\substack{\scriptscriptstyle{p+q=k}}}
\Hom_{\rD^\rb_{\rc}(X)}(\C_{U},\C_{V}[p])\otimes_{\Bbb C}\Hom_{\rD^\rb_{\coh}(\sho_S)}(F, G[q])\ar[r]^(0.55)*[@]{\cong}
&\Hom_{\rD^\rb_{\rc}(\pOS)}(\shf,\shg[k])\ar[u]_*[@]{\cong}\\
\bigoplus\limits_{\substack{\scriptscriptstyle{p+q=k}}}
\Hom_{\rD^\rb(\Mod_{\rc}(X))}(\C_{U},\C_{V}[p])
\otimes_{\Bbb C}\Hom_{\rD^\rb(\coh(\sho_S))}(F, G[q])\ar[r]
\ar[u]_(0.40)*[@]{\cong}
&\Hom_{\rD^\rb(\Mod_{\rc}(\pOS))}(\shf, \shg[k])\ar[u]\\
}
\]}}

Thus we conclude from the bottom diagram that the bottom right vertical arrow is surjective which implies that the canonical map \eqref{FF} is surjective.

Let us now prove the essential surjectivity: the essential image $\shd'$ of $\shi$ is a strictly full triangulated subcategory of $\rD^\rb_{\rc}(\pOS)$
containing $\she$. Thus $\shd'=\langle \she \rangle= \rD^\rb_{\rc}(\pOS)$ if and only if it is thick
which is equivalent to prove that any idempotent in $\shd'$ splits. 
Given $e$ an idempotent in $\shd'$, since $\shi$ is fully faithful, there exists a unique $\epsilon$ idempotent in
$\rD^\rb(\Mod_{\rc}(\pOS))$ such that $\shi(\epsilon)=e$.
According to \cite[Theorem]{LC} $\epsilon$ splits in $\rD^\rb(\Mod_{\rc}(\pOS))$ which implies that
$e$ splits in the essential image of $\shi$.
\end{proof}

\begin{proposition}\label{PN3}
 Let $S$ be
a complex manifold. Assume that there exists a real analytic manifold X
such that the  natural  functor
$$
\shi:\rD^\rb(\Mod_{\rc}(\pOS))\to \rD^\rb_{\rc}(\pOS)$$
 is an equivalence of categories.
Then the functor $E$ in \eqref{E200}
is an equivalence of categories.

\end{proposition}
\begin{proof}
 If $X$ is empty the result is trivial. Let us assume that $X$ is non empty and let us choose  a point $x_0\in X$.
Let 
 $\Mod_{\rc,x_0}(\pOS)$ (resp. $\Mod_{x_0}(\pOS)$) be the full abelian subcategory
of $\Mod_{\rc}(\pOS)$ (resp. $\Mod(\pOS)$) whose objects have support contained in $\{x_0\}\times S$.
In view of \cite[Prop.\,1.7.11]{KS1},
$\rD^\rb(\Mod_{\rc,x_0}(\pOS))$  is equivalent to the full subcategory of
 $\rD^\rb(\Mod_{\rc}(\pOS))$ (resp.  $\rD^\rb(\Mod_{x_0}(\pOS))$ is equivalent to the full subcategory of
 $\rD^\rb(\Mod(\pOS))$) whose objects have cohomology supported on $\{{x_0}\}\times S$.
 Let $\rD^\rb_{\rc,x_0}(\pOS)$ denote the full subcategory of $\rD^\rb_{\rc}(\pOS)$
whose objects have support contained in $\{x_0\}\times S$, by the previous remark this is equivalent to 
$\rD^\rb_{\rc}(\Mod_{x_0}(\pOS))\simeq \rD^\rb_\coh(\sho_S)$.
Then again $\shi$ induces an equivalence
$\shi_{x_0}:\rD^\rb(\Mod_{\rc,x_0}(\pOS))\to \rD^\rb_{\rc,x_0}(\pOS).$ 
We have the following commutative diagram where the vertical arrows are given by the fully faithful functor $\C_{\{x_0\}}\boxtimes \cbbullet $. 
$$
\xymatrix@R=10pt{
\rD^\rb(\coh(\sho_S))\ar[r]\ar[d]  &\rD^\rb_{\coh}(\sho_S)\ar[d] \\
\rD^\rb(\Mod_{\rc,x_0}(\pOS))\ar[r]^(0.6){\shi_{x_0}} & \rD^\rb_{\rc,x_0}(\pOS)
}
$$
The vertical and the bottom arrows are equivalences thus the top arrow is also an equivalence of categories.
\end{proof}
\subsection{The algebraic view point}Let $X$ be a complex manifold and let us assume that $S$ is the analytification of a smooth  complex algebraic variety $S^{\mathrm{alg}}$.
We can endow the product $X\times S^{\mathrm{alg}}$ with a structure of  ringed space. Let $\pi_X:X\times S^{\mathrm{alg}}\to S^{\mathrm{alg}}$ be the projection and let $h: X\times S \to X\times S^{\mathrm{alg}}$
be the canonical analytification morphism of  ringed spaces.
This analytification procedure has been used in \cite{Wu21}.

We can easily translate the notion of
$S\textup{-}\R\textup{-}$constructible sheaves introduced in \S\ref{pre}
replacing open neighborhoods on $S$ with Zariski open neighborhoods on $S^\alg$ and $\sho_S\textup{-}$coherent modules with
$\sho_{S^\alg}\textup{-}$coherent ones. We
denote by $\Mod_{\rc}(\pi_X^{-1}\sho_{S^\alg})$ the abelian category of
$S^\alg\textup{-}\R\textup{-}$constructible sheaves of $\pi_X^{-1}\sho_{S^\alg}$-modules.

We remark that Lemma~\ref{L1} holds true if we replace
$\she$ with $\she^\alg$ the family of locally finite direct sums of objects  $\C_{U}\boxtimes L$ where $U$ is relatively compact contractible subanalytic open in $X$ and $L$ is a 
locally free $\sho_{S^\alg}$-module of finite rank.
Then $\she^\alg$ classically generates $\rD^\rb_{\rc}(\pi_X^{-1}\sho_{S^\alg})$. 

The following result is an easy consequence of Lemma~\ref{L1} and Theorem~\ref{Cor1}.

\begin{corollary}
If $S^\alg$ is projective
$$
h:\Mod_{\rc}(\pi_X^{-1}\sho_{S^\alg})\to \Mod_{\rc}(\pOS)$$
 is an equivalence of categories which induces the following equivalences of categories:
 $$
\rD^\rb_{\rc}(\pi_X^{-1}\sho_{S^\alg})\to \rD^\rb_{\rc}(\pOS)
 $$
$$
\rD^\rb_{\cc}(\pi_X^{-1}\sho_{S^\alg})\to \rD^\rb_{\cc}(\pOS).
 $$

\end{corollary}
We think that the preceding results would be the starting point for a an algebraic relative Riemann-Hilbert equivalence assuming $S^\alg$ projective.

\section{Relative Integral transforms}\label{IT}

This section only depends on the results in 1.a.
 We will adapt to the relative case, without insisting in details, the theory of integral transforms with kernels for sheaves and $\shd$-modules, essentially adapting the formalisms introduced in \cite{KS1} and in \cite{KS2}. 
For instance, the notions of composition of relative kernels for sheaves or for relative $\shd$-modules is a natural adaptation of the corresponding ones in the absolute case (see \cite[Defs. B1, B4]{D-S}). 
In accordance with the notations in \cite{FMFS2}, for a morphism $\phi=(f,g): Y\times T\to X\times S$ we set $L{\phi} ^*:=f^{-1}Lg^*$ where  $f^{-1}$ stands for $(f,\id)^{-1}$ and $g$ stands for $(\Id, g)$ while the notion of $R\phi_*$ is clear.

Following the classic construction of integral transforms, let us consider real analytic manifolds $X$, $Y$ and complex manifolds $S$, $T$ together with the relative framework concerned by $p_X:X\times S\to S$, $p_Y:Y\times T\to T$ and $p_{X\times Y}: X\times Y\times S\times T\to S\times T$. 
We denote by $q_X$ the projection $X\times Y\times S\times T\to X\times S$ and by  $q_Y$ the projection $ X\times Y\times S\times T\to Y\times T$. 
Thus $q_X=q'_X\times \pi$, $q_Y=q'_Y\times \pi'$ where $q'_X$ is the projection $X\times Y\to X$, $\pi$ is the projection $S\times T\to S$,  $q'_Y$ is the projection $X\times Y\to Y$ and $\pi'$ is the projection $S\times T\to T$. 
\begin{definition}\label{DRIT}Given $K\in\rD^\rb_{\rc}(p^{-1}_{X\times Y}\sho_{S\times T})$ such that $q_Y$ is proper in $\supp K$, then, in view of Propositions \ref{PN1} and \ref{PN2}, the functor \begin{equation}\label{E111}\Phi_K: F\mapsto Rq_{Y*}(K\otimes^L_{p^{-1}_{X\times Y}\sho_{S\times T}}q^*_XF)
\end{equation}
 is well defined from $\rD^\rb_{\rc}(\pOS)$ to $\rD^\rb_{\rc}(p^{-1}_Y\sho_{T})$ and we call it the relative integral transform with kernel $K$.
 \end{definition}
Clearly, if $X$ and $Y$ are complex manifolds and $K$ belongs to $\rD^\rb_{\cc}(p^{-1}_{X\times Y}\sho_{S\times T})$ then $\Phi_{K}$ restricts as a functor from 
$\rD^\rb_{\cc}(p_X^{-1}\sho_S)\to \rD^\rb_{\cc}(p_Y^{-1}\sho_T)$.

\subsection{Example with $S$ fixed: relative projective duality}\label{Sfix}

Recall that, for $F\in\rD^\rb(\C_X)$ we note $$\bD'(F):=\Rhom(F, \C_X).$$
Let $X$ and $Y$ be real analytic manifolds and let us suppose that $S=T$. We identify $S$ with the diagonal $\Delta_S\subset S\times S$.

As an example of \eqref{E111} we find the relative Radon transform (with $S=T$ and $K=\C_{Z\times S}\boxtimes\sho_S$ where $Z$ is a relatively compact subanalytic subset of $X\times Y$). 
This study is completely similar to the absolute case treated for instance in \cite{S2}, \cite{D-S}. For the sake of applications we will only treat here the relative versions of real and complex projective duality (\cite{D-S}). For such a purpose we need to recall some notations and results in loc.cit..

\textit{Relative real projective duality.} Let $\rm{P}$ be a $n$-dimensional real projective space and let $\rm{P}^*$ be its dual projective space. We denote by $\rm{A}$ the real analytic incidence hypersurface of $\rm{P}\times \rm{P}^*$ defined by $\rm{A}=\{(x,\xi); \langle x,\xi\rangle=0\}$ and set $\Omega=\rm{P}\times\rm{P}^*\setminus \rm{A}$. We note $r: \rm{P}\times\rm{P}^*\to \rm{P}^*\times\rm{P}$ the morphism $r(x,y)=(y,x)$. Then, as proved in \cite[Lem.~2.3]{CM}, the integral transform with kernel $K:=\C_{\Omega}$ defines an equivalence of categories
from $\rD^\rb_{\rc}(\rm{P})$ to $\rD^\rb_{\rc}(\rm{P}^*)$, a quasi-inverse being the integral transform with kernel $^tK[n]:=\bD'(Rr_*\C_{\Omega})[n]$.
Applying the variant of Nakayama's lemma (\cite[Prop.~2.2]{MFCS1}), we easily derive the following result in the relative setting.

\begin{lemma}\label{RIT2}
The relative integral transform $\Phi_{K}$ where $K=\C_{\Omega}\boxtimes \sho_S$ induces an equivalence of categories from $\rD^\rb_{\rc}(p^{-1}_{\rm{P}}\sho_S)$ to $\rD^\rb_{\rc}(p^{-1}_{\rm{P}^*}\sho_S)$, a quasi-inverse being the relative integral transform $\Phi_{K'}$ with $K'=\bD'(Rr_*\C_{\Omega})[n]\boxtimes \sho_S$.
\end{lemma}

\textit{Relative complex projective duality.} Let us denote by $\mathbb{P}$ a complex $n$-dimensional projective space and by $\mathbb{P}^*$ its dual projective space.
Let $\mathbb{A}$ denote the analytic incidence hypersurface of $\mathbb{P}\times\mathbb{P}^*$ $$\mathbb{A}=\{(x,\xi); \langle x, \xi\rangle=0\}$$ and let now $\Omega:=\mathbb{P}\times\mathbb{P}^*\setminus \mathbb{A}$.
According to \cite[page 39, lines-6, -7]{D-S}, the integral transform with kernel $K:=\C_{\Omega}$ induces an equivalence of categories from $\rD^\rb_{\cc}(\mathbb{P})$ to $\rD^\rb_{\cc}(\mathbb{P}^*)$ 
(proved in \cite{KT})
with quasi-inverse the integral transform with kernel $\bD'(Rr_*\C_{\Omega})[2n]$.
 Similarly to the preceding Lemma, the variant of Nakayama's lemma entails:
\begin{proposition}\label{RIT3}
The relative integral transform $\Phi_K$ with $K=\C_{\Omega}\boxtimes \sho_S$ induces an equivalence of categories from $\rD^\rb_{\cc}(p^{-1}_{\mathbb{P}}\sho_S)$ to $\rD^\rb_{\cc}(p^{-1}_{\mathbb{P}^*}\sho_S)$, a quasi-inverse being the relative integral transform $\Phi_{K'}$ with $K'=\bD'(Rr_*\C_{\Omega})[2n]\boxtimes \sho_S$. 
\end{proposition}

\subsection{Fourier-Mukai transform in the base space}
 In this section we assume that $S$ is a complex abelian variety of dimension g.  
We denote by $\hat{S}$ the dual abelian manifold and by $\mathcal{P}$  the normalized Poincar\'e line bundle on $S\times\hat{S}$ which we identify with the associated locally free of rank one $\sho_{S\times \hat{S}}$-module. 
 
 One also denotes by $q_S$ (resp. by $q_{\hat{S}}$) the projection $S\times \hat{S}\to S$ (resp. the projection $S\times \hat{S}\to \hat{S}$).
 
We note that Fourier-Mukai transforms are defined in more general situations. Here we confine to the coherent case sufficient for our applications. The notations in this section follow \cite{Muk81}, being slightly different from those used in the Introduction.
 
 \begin{itemize} \item{One denotes by $\rm{RS}$ the Fourier-Mukai integral transform (with kernel $\shp$) from $\rD^\rb_{\coh}(\sho_{\hat{S}})$ with values in $\rD^\rb_{\coh} (\sho_{S})$ given by 
 $$\hat{F}\mapsto Rq_{S*}(\shp\otimes_{\sho_{S\times\hat{S}}}q^*_{\hat{S}}\hat{F})\simeq Rq_{S*}(\shp\otimes_{q^{-1}\sho_{\hat{S}}}q^{-1}\hat{F}).$$ 
  }
  \item{One denotes by $\rm{R\hat{S}}$ the Fourier-Mukai integral transform (with kernel $\shp$) from $\rD^\rb_{\coh}(\sho_{S})$ with values in $\rD^\rb_{\coh} (\sho_{\hat{S}})$ given by 
 $${F}\mapsto Rq_{\hat{S}*}(\shp\otimes_{\sho_{S\times\hat{S}}}q^*_{S}F)\simeq Rq_{\hat{S}*}(\shp\otimes_{q^{-1}\sho_{S}}q^{-1}F).$$}
 \item{We have canonical isomorphisms of functors (inversion formulas) proved by H. Liu in \cite[Th. 4.1.1]{HL}
 \begin{equation}\label{EHL1}\rm{RS}\circ \rm{R\hat{S}}\simeq (-\Id_S)^*[-g]
 \end{equation}
 \begin{equation}\label{EHL2}\rm{R\hat{S}}\circ \rm{R S}\simeq (-\id_{\hat{S}})^* [-g]
 \end{equation} thus $\rm{R S}$ and $\rm{R\hat{S}}$ are equivalences of categories.
 }

 \end{itemize}
 \subsubsection{Relative Fourier-Mukai transform} Let $X$ be a real analytic manifold. We shall keep the notation $q_{S}$ (resp. $q_{\hat{S}}$) for the projection $X\times S\times\hat{S}\to X\times S$ (resp. the projection $X\times S\times\hat{S}\to X\times \hat{S}$). In view of Propositions \ref{PN1} and \ref{PN2} the following integral transformations (with kernel $p^{-1}_{X}\shp$) are well defined: 
 
 $\rm{RS_X}: \rD^\rb_{\rc}(p^{-1}_{X}\sho_{\hat{S}})\to \rD^\rb_{\rc}(p^{-1}_{X}\sho_S)$ (resp. $\rm{R\hat{S}_X}: \rD^\rb_{\rc}(p^{-1}_{X}\sho_S)\to \rD^\rb_{\rc}(p^{-1}_{X}\sho_{\hat{S}})$) given by 
 \begin{equation}\label{EHL3}\hat{F}\mapsto Rq_{S*}(p^{-1}_{X}\shp\otimes_{p^{-1}_{X}\sho_{S\times\hat{S}}}q^*_{\hat{S}}\hat{F})
 \end{equation} (resp. by \begin{equation}\label{EHL4}F\mapsto Rq_{\hat{S}*}(p^{-1}_{X}\shp\otimes_{p^{-1}_{X}\sho_{S\times\hat{S}}}q^*_{S}F)
 \end{equation}
 Because the kernel of both transformations is $p_X^{-1}\mathcal{P}$, Theorem 2.2 in \cite{Muk81} has a  straightforward translation in the relative setting:
 \begin{proposition}\label{PLH1}
The inversion formulas \eqref{EHL1} and \eqref{EHL2} hold true for the transformations \eqref{EHL3} and \eqref{EHL4}, in particular the latter define equivalences of categories. 
\end{proposition}

 \begin{remark}\label{RRIT6} \textit{The relative $\C$-constructible case.}
In view of Propositions \ref{PN1} and \ref{PN2}, the transformations \eqref{EHL3} and \eqref{EHL4} restrict to $\rD^\rb_{\cc}(\pOS)$ and $\rD^\rb_{\cc}(p^{-1}_X\sho_{\hat{S}})$. Thus the inversion formulas \eqref{EHL1} and \eqref{EHL2}  also induce  equivalences between these categories. 
  \end{remark}
  
\subsection{Integral transforms for relative regular holonomic $\shd$-modules }
Let us consider the situation described in the beginning of this section.
We assume that $S$ is a complex projective algebraic variety thus the projection $S\times T\to T$ is a projective morphism.

Let $\shh$ be a complex of $\shd_{X\times Y\times S\times T/S\times T}$ with 
and regular holonomic cohomology and assume that the projection $q_Y: X\times Y\times S\times T\to Y\times T$ is proper on $\supp \shh$.

Then, according to the stability of regular holonomicity under pullback and proper pushforward (\cite[Th. 1]{FMFS2} and \cite[Cor. 2.4]{MFCS2}), pullback  and projective pushforward in the base space (\cite[Prop. 3.10, Prop. 3.15]{FMFS2}) and by \cite[Prop. 4.34]{FMFS2}, the correspondence 
\begin{equation}\label{RITDM}
\Phi_{\shh}:\shm\mapsto \Dq_{Y*}(\shh\otimes^L_{\sho_{X\times Y\times S\times T}}q_X^{*}\shm)
\end{equation}
is a well defined functor from $\rD^\rb_{\rhol, \rm{good}}(\DXS)$ to $\rD^\rb_{\rhol,   \rm{good}}(\shd_{Y\times T/T})$. 
\begin{definition}\label{DRIT2}
The functor \eqref{RITDM} is called the relative integral transform with kernel $\shh$. \end{definition}

\begin{example}\label{ERIT}
Let $A$ be an abelian variety and let $A^{\sharp}$ be the moduli space of line bundles with integrable connection on A. The generalized Fourier-Mukai integral transform $\rm{FM}$ (\cf Introduction)  is a particular case of  \eqref{RITDM} with $S=Y=\{pt\}$ ($S$ is trivially projective), $X=A$, $T=A^{\sharp}$. The results in \cite{FMFS2} explain why this integral transform is well defined with the simpler assumption of coherence, that is for $\shm\in \rD^\rb_{\coh}(\shd_A)$. The main point is that its kernel $\mathcal{P}$ is a flat relative connection so \cite[Prop. 2.11 (ii)]{FMFS2} implies that $\rm{FM}(\shm)$ is an object in $\rD^\rb_{\coh}(\shd_{\{pt\}\times A^{\sharp}/A^{\sharp}})=\rD^\rb_{\coh}(\sho_{A^{\sharp}})$. We remark that by \cite[Th.~18.1]{Sch15} the functor $\rm{FM}$ is t-exact with respect to the standard t-structure on $\rD^\rb_{\coh}(\shd_A)$
and the m-t-structure on $\rD^\rb_{\coh}(\sho_{A^{\sharp}})$ (see \cite[Lem.~17.4]{Sch15}). 
This result is connected to \cite[Th.~4.1]{FMF1}.
\end{example}

\remark 
Of course, if in Definition \ref{DRIT2}, $S=T$ and $\shh$ is supported on $X\times Y\times\Delta_S$ (see section \ref{Sfix}), $S$ does not need to be projective since $q_Y|_{\supp \shh\otimes^L_{\sho_{X\times Y\times S\times T}}q_X^*\shm}$ is of the form $q'_Y\times \Id_S$.

Coming back to the relative complex projective duality (\cf section \ref{Sfix}), by the relative Riemann-Hilbert equivalence and  Proposition~\ref{RIT3}, we conclude an equivalence of categories 
\begin{equation}\label{RIT4}\rD^\rb_{\rhol}(\shd_{\mathbb{P}\times S/S})\to \rD^\rb_{\rhol}(\shd_{\mathbb{P}^*\times S/S})
\end{equation}
One easily proves that the latter when restricted to $\rD^\rb_{\rhol,\rm{good}}(\shd_{\mathbb{P}\times S/S})$ is the integral transform $\Phi_{\shh}$ with $$\shh=\RH^S_{\mathbb{P}\times\mathbb{P}^*}(\C_{\Omega}\boxtimes\sho_S)[-2n]\simeq \tho(\C_{\Omega}\times S, \sho_{\mathbb{P}\times\mathbb{P}^*\times S})$$ The main step is to apply the variant of Nakayma's lemma reducing to the absolute case which was proved in \cite{D-S}.

\begin{example} \textit{Relative Fourier-Mukai transform.}\label{ERIT2} According to the relative Riemann-Hilbert correspondence, in view of Remark \ref{RRIT6}, transformations \eqref{EHL3}, \eqref{EHL4} induce equivalences of the categories  $\rD^\rb_{\rhol}(\DXS)$ and $\rD^\rb_{\rhol}(\shd_{X\times\hat{S}/\hat{S}})$ which we denote by $\rm{RS_X}$ and $\rm{R\hat{S}_X}$ respectively. We will study $\rm{R\hat{S}}_X$ as an equivalence of categories 
\begin{equation}\label{RIT7}
\rD^\rb_{\rhol,\rm{good}}(\DXS)\to \rD^\rb_{\rhol,\rm{good}}(\shd_{X\times\hat{S}/\hat{S}})
\end{equation}
We shall prove that \eqref{RIT7} is the integral transform $\Phi_{\shl}$ with 
$\shl=
\RH^{S\times \hat{S}}_{X}(p_X^{-1}(\mathcal{P}\otimes_{\sho_{S\times\hat{S}}}\Omega^{-1}_{S\times \hat{S}/\hat{S}}[-d_S]))[-d_X]$.
We start by noting that, since $\mathcal{P}\otimes_{\sho_{S\times\hat{S}}}\Omega^{-1}_{S\times \hat{S}/\hat{S}}$ is a locally free  $\sho_{S\times\hat{S}}$-module of finite rank (and hence coherent), we have (\cf\cite[Cor. 3.24]{MFCS2}) 
$\RH^{S\times \hat{S}}_{X}(p_X^{-1}(\mathcal{P}\otimes_{\sho_{S\times\hat{S}}}\Omega^{-1}_{S\times \hat{S}/\hat{S}}[-d_S]))[-d_X]\simeq \Rhom_{p_{X}^{-1}\sho_{S\times \hat{S}}}(
p_X^{-1}(\mathcal{P}\otimes_{\sho_{S\times\hat{S}}}\Omega^{-1}_{S\times \hat{S}/\hat{S}}[-d_S]), \sho_{X\times S\times\hat{S}})$. By the same reason, $\shl$ is a relative flat holomorphic connection therefore it is a good $\shd_{X\times S\times\hat{S}/S\times \hat{S}}$-module. 

We aim to check that, for any $\shm\in\rD^\rb_{\rhol, \rm{good}}(\DXS)$, we have a functorial isomorphism in $\rD^\rb_{\rhol, \rm{good}}(\DXS)$.
\begin{equation}\label{ERFM2}
\Phi_{\shl}\shm\simeq
\RH_X^{\hat{S}}(\rm{R\hat{S}}_X(\pSol \shm))\end{equation}
Let us set $\pSol \shm=F$ which is equivalent to $\shm\simeq \RH_X^S(F)$.
The right side of \eqref{ERFM2} is $$\RH_X^{\hat{S}}(Rq_{\hat{S}*}(p_X^{-1}\mathcal{P}\otimes^L_{p_X^{-1}\sho_{S\times \hat{S}}}q^*_SF))$$
which is isomorphic to 
$Rq_{\hat{S}*}\RH_X^{S\times\hat{S}}(
p^{-1}_{X}(\mathcal{P}\otimes_{\sho_{S\times\hat{S}}}\Omega^{-1}_{S\times \hat{S}/\hat{S}}[-d_S])
\otimes^L_{p_X^{-1}\sho_{S\times\hat{S}}}q_S^*F)$
by Proposition~\ref {Ldirectim}.

The left side of \eqref{ERFM2} is 
\begin{align*}
&
Rq_{\hat{S}*}(\RH^{S\times \hat{S}}_{X}(
p_X^{-1}(\mathcal{P}\otimes_{\sho_{S\times\hat{S}}}\Omega^{-1}_{S\times \hat{S}/\hat{S}}[-d_S]))[-d_X]
\otimes^L_{\sho_{X\times S\times\hat{S}}}q_S^*\RH^S_X(F))\\
&\simeq
Rq_{\hat{S}*}(\Rhom_{p^{-1}\sho_{S\times\hat{S}}}(p^{-1}_{X}(\mathcal{P}\otimes_{\sho_{S\times\hat{S}}}\Omega^{-1}_{S\times \hat{S}/\hat{S}}[-d_S]), \sho_{X\times S\times\hat{S}})\otimes^L_{\sho_{X\times S\times\hat{S}}}q_S^*\RH^S_X(F))
\\
%\end{align*}
%\begin{align*}
%&Rq_{\hat{S}*}(\Rhom_{p^{-1}\sho_{S\times\hat{S}}}(p^{-1}_{X}(\mathcal{P}, \sho_{X\times S\times\hat{S}})\otimes^L_{\sho_{X\times S\times\hat{S}}}q_S^*\RH^S_X(F))\\
%&\simeq Rq_{\hat{S}*}(\Rhom_{p^{-1}\sho_{S\times\hat{S}}}(p^{-1}_{X}\mathcal{P}, q_S^*\RH^S_X(F)))\\
& \simeq Rq_{\hat{S}*}\RH_X^{S\times\hat{S}}(
p^{-1}_{X}(\mathcal{P}\otimes_{\sho_{S\times\hat{S}}}\Omega^{-1}_{S\times \hat{S}/\hat{S}}[-d_S])
\otimes^L_{p_X^{-1}\sho_{S\times\hat{S}}}q_S^*F)\\
\end{align*} according to Lemma \ref{LRFM}, Proposition \ref{Linvim} and by the coherence of $\mathcal{P}$. The assumption on $\shm$ implies that $q_S^*\shm\simeq q_S^*\RH_X^S(F)$ is also a good $\shd_{X\times S\times\hat{S}/S\times \hat{S}}$-module thus, by the  previous argument, $\RH_X^{S\times\hat{S}}(p^{-1}_{X}\mathcal{P}\otimes^L_{p_X^{-1}\sho_{S\times\hat{S}}}q_S^*F)$ has $q_{\hat{S}}$-good cohomology. %The statement then follows by Proposition~\ref {Ldirectim}.
\end{example}

\begin{remark}\label{R18}When $X=\{pt\}$, Example \ref{ERIT2} is just the classical Fourier-Mukai transform for coherent $\sho$-modules up to replacing the kernel $\mathcal{P}$ by its dual $\bD(\mathcal{P}):=\Rhom_{\sho_{S\times \hat{S}}}(\mathcal{P},\sho_{S\times \hat{S}})$.
\end{remark}
 
 \begin{example}\label{ERRIT8}Let us consider the following diagram of equivalences of categories $$
\xymatrix@R=10pt{
\rD^\rb_{\cc}(p^{-1}_{\mathbb{P}}\sho_S)\ar[r]^{\Phi_{\shk}} &\rD^\rb_{\cc}(p^{-1}_{\mathbb{P}^*}\sho_S)\ar[r]^{\rm{R\hat{S}}_{\mathbb{P}^*}}&\rD^\rb_{\cc}(p^{-1}_{\mathbb{P}^*}\sho_{\hat{S}})\ar[d]^{\RH^{\hat{S}}_{\mathbb{P}^*}} \\
\rD^\rb_{\rhol}(\shd_{\mathbb{P}\times S/S})\ar[r]\ar[u]^{\pSol_{\mathbb{P}}} & \rD^\rb_{\rhol}(\shd_{\mathbb{P}^*
\times S/S})\ar[r]\ar[u]^{\pSol_{\mathbb{P}^*}}&\rD^\rb_{\rhol}(\shd_{\mathbb{P}^*\times \hat{S}/\hat{S}})
}
$$

Since by construction the bottom arrows are equivalences we conclude that 
the composition $\rm{R\hat{S}}_{\mathbb{P}^*}\circ \Phi_{\shk}$ induces an equivalence of categories
$$\rD^\rb_{\rhol}(\shd_{\mathbb{P}\times S/S})\to \rD^\rb_{\rhol}(\shd_{\mathbb{P}^*\times \hat{S}/\hat{S}})$$ which is nothing more then the integral transform with kernel $\shl\circ\shh$ when restricted to $\rD^\rb_{\rhol, \rm{good}}(\shd_{\mathbb{P}\times S/S})$.
\end{example}
\section{Relative constructible functions}\label{SCF}

After the introduction and study of the category of relative $\R$-constructible sheaves along several papers (\cite{MFCS1}, \cite{MFCS2}, \cite{FMFS1}, \cite{FMFS2}) the purpose of this section is to answer to the natural question of existence of a corresponding notion of relative constructible functions together with its essential properties.

\subsection{$\sha$-constructible functions}
Let $\sha$ be a commutative unital ring.
\begin{definition}\label{DfSc}
A function $\varphi: X\to \sha$ is called $\sha$-constructible if:
\begin{enumerate}
\item{For each $y\in \sha$, $\varphi^{-1}(y)$ is subanalytic.}
\item{The family $(\varphi^{-1}(y))_{y\in \sha}$ is locally finite.}
\end{enumerate}
\end{definition}
In particular, given $\varphi\in \CF(X)$ and $y\in\sha$, $$(y\cdot\varphi)(x):=y\cdot\varphi(x)$$
is a well defined $\sha$-constructible function.  
\begin{remark}
In the case of $\sha=\Z$, the previous definition reduces to that of constructible function in \cite[(9.7.1)]{KS1}. We recall that the ring of constructible functions on $X$ is denoted by $\CF(X)$ in \cite{KS1}.
\end{remark}
 As a consequence, the set of $\sha$-constructible functions is naturally endowed with a structure of an $\sha$-algebra which we denote by $\CF^{\sha}(X)$. The presheaf $$U\mapsto \CF^{\sha}(U)$$ defines a sheaf which we denote by 
  $\mathcal{C}\mathcal{F}_X^{\sha}$ (or $\mathcal{C}\mathcal{F}^{\sha}$ if there is no ambiguity) which is soft.
In this way we can identify $\CF(X)$
with a subring of $\CF^{\sha}(X)$ by associating to $\varphi \in \CF(X)$ the $\sha$-constructible function $1\cdot\varphi$.
We can adapt the contents of \cite[Pag.\,399]{KS1} as follows:
\begin{proposition}\label{P1fSc}
A function $\varphi$ is $\sha$-constructible if and only if there exists a $\mu$-stratification $X=\bigsqcup_{\alpha} X_{\alpha}$ such that $\varphi|_{X_{\alpha}}$ is constant for each $\alpha$.
\end{proposition}
\begin{proof}
If $\varphi$ is $
\sha$-constructible then, according to \cite[Th.\,8.3.20]{KS1}, we can choose a $\mu$-stratification $X=\bigsqcup_{\alpha} X_{\alpha}$  refining the locally finite family of subanalytic subsets $(\varphi^{-1}(y))_{y\in \sha}$. Thus for each $\alpha$, $\varphi|_{X_{\alpha}}$ is constant.
The converse is clear.
\end{proof}
Another characterization is the following. We denote by ${\bf{I}}_Y$ the characteristic function of a subset $Y$.
\begin{proposition}\label{P2fSc}
The function $\varphi$ is $\sha$-constructible if and only if there exists a locally finite covering $X=\bigcup_{\beta} Y_{\beta}$ where each $Y_{\beta}$ is compact subanalytic and contractible, and a family $(y_{\beta})_{y_{\beta}\in\sha}$ such that 
\begin{equation}\label{ast}
\varphi=\sum_{\beta} y_{\beta}\cdot{\bf{I}}_{Y_{\beta}}
\end{equation}
\end{proposition}
\begin{proof}
If $\varphi$ is $\sha$-constructible, considering a $\mu$-stratification as in Proposition~\ref{P1fSc}, we apply the Triangulation Theorem (\cf \cite[Th.\,8.3.25]{KS1}) to obtain a subanalytic triangulation of $X$ such that each $X_{\alpha}$ is a union of simplices $Y_{\alpha_j}$ (hence compact subanalytic and contractible). We define
$y_{\alpha_j}=\varphi(x)$ with $x\in Y_{\alpha_j}$.
Thus $Y_{\alpha_j}$  with $\beta=\alpha_j$ and $y_{\beta}=y_{\alpha_j}$
satisfy the desired conditions.
Conversely, given the family $(Y_{\beta})$ satisfying \eqref{ast}, we can apply again \cite[Th.\,8.3.20]{KS1} to obtain, by refining $(Y_{\beta})$, a $\mu$-stratification in the conditions of Proposition~\ref{P1fSc}.
\end{proof}
Here and in the sequel, unless ambiguity, $[\cbbullet]$ will denote the class of an object $\cbbullet$ in a given Grothendieck group.

\begin{definition}[Index over $\sha$]\label{D5base}
An index over $\sha$ is a ring morphism
$$\psi^{\sha}(\cbbullet): \rK_0(\rD^\rb_{\rc}(\pOS))\to \CF^{\sha}(X).$$
\end{definition}

\subsection{$S$-constructible functions}\label{S}
For the sake of simplicity, we shall denote by $\rK_0(S)$ the Grothendieck group of  the triangulated category $\rD^\rb_\coh(\sho_S)$ and also denote by ${\rK}_{\rc}(\pOS)$ the Grothendieck group $\rK_0(\rD^\rb_{\rc}(\pOS))$.
We notice that, by \cite[Prop. A.9.5]{Ac}, the natural morphism $\rK_0(\coh(\sho_S))\to \rK_0(S)$ (resp. $\rK_0(\Mod_{\rc}(\pOS))\to \rK_{\rc}(\pOS)$) is an isomorphism of groups. In both cases the inverse morphism is given by 
\begin{equation}\label{Eqiso2}
\begin{matrix}
\hfill [F]&\longmapsto & \sum_j(-1)^j[\shh^jF].\\
\end{matrix}
\end{equation}

We recall that 
$\rK_0(S)$ (resp. $\rK_{\rc}(\pOS)$) is a commutative ring whose
product is defined by: 
\begin{equation}\label{EN}
[F]\cdot [G]:=[F\otimes^L_{\sho_S}G] 
\end{equation}
(resp. \begin{equation}\label{EN2a}
[F]\cdot [G]:=[F\otimes^L_{\pOS}G] )
\end{equation}
for $F,\,G\in \rD^\rb_\coh(\sho_S)$ (resp. $F, G\in \rD^\rb_{\rc}(\pOS)$) .
 Clearly, $[\sho_S]$ (resp $[\pOS]$) is the unit. 
In the sequel we replace the commutative algebra $\sha$ by $\rK_0(S)$ and, for simplicity, we also replace the notations $ \CF^{\sha}$ by $ \CF^S$ and $\sha$-constructible by $S$-constructible. Thus from Definition~\ref{DfSc} we obtain the notion of $S$-constructible function (or relative constructible function).

We define the morphism $\chi^S: \rK_{\rc}(\pOS)\to \CF^S(X)$
\begin{equation}\label{E1}\chi^S([F])(x)=[F|_{\{x\}\times S}]_{\rK_0(S)}=\sum_i (-1)^i [\shh^i(F|_{\{x\}\times S})]_{\rK_0(S)}
\end{equation} for any $[F]\in{\rK}_{\rc}(\pOS)$. This definition is clearly independent of the choice of the representative of $[F]$ in $\rD^\rb_{\rc}(\pOS)$.

\begin{definition}\label{Def:EPI}
For $[F]\in {\rK}_{\rc}(\pOS)$, $\chi^S(F)$ is the Euler-Poincar\'e index of $[F]$.
\end{definition}
 
\begin{example}\label{EX22}
Let us recall \cite[Ex. 3.12]{MFCS1}.
Let $X$ be the open unit disc in $\C$ with coordinate $x$ and
let $S$ be a connected open set of $\C$ with coordinate $s$. 
Let $\phi : S\to \C$
be a non constant holomorphic function on $S$ and consider the holonomic
$\DXS$-module $\shm = \DXS/\DXS P$, 
with $P(x,\partial_ x, s) = x\partial_x-\phi(s)$. 
Then $\shm$ admits the resolution
$$0\to \DXS\underset{\cdot P}{\to}\DXS\to \shm\to 0$$ The complex $\Sol M$ is
represented by $0\to \sho_{X\times S}\overset{P.}{\to}
\sho_{X\times S}\to 0$ (terms in degrees 0 and 1).
Consider the stratification $X_1 = X\setminus \{0\}$ and $X_0 = \{0\}$ of $X$. Then
$H ^{0}\Sol \shm|_{X_1\times S}$ is a locally constant sheaf of free $\pOS$-modules generated
by a local determination of $x^{\phi(s)}$, and $H^1 \Sol \shm|_{X_1\times S} = 0$. Also $H^{0} \Sol \shm|_{X_0\times S} = 0$ and $H^1 \Sol \shm|_{X_0\times S}$ is a skyscraper sheaf on $X_0\times S$ supported on $\{s\in S\colon\phi(s)\in \Z \}$, a torsion module. Then $\chi^S(\Sol \shm)(0)=0$
 and $\chi^S(\Sol\shm)(x)=[\sho_S]=1,\,x\neq 0$.
\end{example}

\begin{remark}
 If $\dim S=0$ the previous definition reduces to the definition of the  Euler-Poincar\'e index $\chi:\rK_{\rc}(X)\to \CF(X)$. 
\end{remark}

\begin{remark}\label{Rem:phiG}
In particular given $\varphi=\chi([F])$ with $[F]\in \rK_\rc(X)$ and $G\in\rD^\rb_\coh(\sho_S)$ then
$\varphi[G]=\chi^S[q_X^{-1}F\otimes p^{-1}_XG]$. Here $q_X$ denotes the projection of $X\times S$ on $X$.
\end{remark}

\begin{remark}\label{R1}
Let us recall that any triangulated functor between two triangulated categories induces a functor between their $\rK_0$.
Let $Z$ be a subanalytic subset of $X$ and let $G$, $H$ be complexes of coherent $\sho_S$-modules.
Then,  $[G]_{\rK_0(S)}=[H]_{\rK_0(S)}$ implies $[\C_Z\boxtimes G]_{{\rK}_{\rc}(\pOS)}=[\C_Z\boxtimes H]_{{\rK}_{\rc}(\pOS)}$. 
Similarly, if $[\C_Z]_{\rK_{\rc}(X)}=[\C_{Z'}]_{\rK_{\rc}(X)}$ for some subanalytic locally closed subsets $Z$ and $Z'$ in $X$, then $[\C_Z\boxtimes G]_{{\rK}_{\rc}(\pOS)}=[\C_{Z'}\boxtimes G]_{{\rK}_{\rc}(\pOS)}$.
We may of course replace $\C$ by any complex of finite dimensional vector spaces.
 \end{remark}

We have the variant of \cite[Th.\,9.7.1]{KS1}:
\begin{theorem}\label{T1fSc}
Let $X$ be a real analytic manifold and let $S$ be a complex manifold. 
Then the morphism $\chi^S$ in \eqref{E1} is an isomorphism. \end{theorem}
\begin{proof}
Isomorphism \eqref{Eqiso2} shows that, for $F\in \rD^\rb_{\rc}(\pOS)$, $[F]=0$ if and only if $$\sum_{j\,odd}[\shh^jF]=\sum_{j\,even}[\shh^jF].$$
Similarly to the proof in \cite[Th.\,9.7.1]{KS1}, any $u\in \rK_{\rc}(\pOS)$ may be represented by a single complex $F\in \rD^\rb_{\rc}(\pOS)$ (i.e. $u=[F]$). 
We shall follow the line of the proof of \cite[Th.\,9.7.1]{KS1} in the absolute case  from which we take the notations.

(a) Let us prove that $\chi^S$ is surjective. Let $\varphi\in \CF^S(X)$. According to Proposition \ref{P2fSc} we may choose a subanalytic stratification $X=\bigsqcup_{\alpha\in A} X_{\alpha}$ and coherent complexes $G_{\alpha}\in \rD^\rb_{\coh}(\sho_S)$ such that $\varphi=\sum_{\alpha\in A}{\bf{I}}_{X_{\alpha}}[G_{\alpha}]$. 
Then, $F=\bigoplus_{\alpha\in A}\C_{X_{\alpha}}\boxtimes G_{\alpha}$ satisfies $\chi^S([F])=\varphi$.

(b) Let us now prove the injectivity of $\chi^S$. Let $u\in \rK_{\rc}(\pOS)$. As in the proof 
of \cite[Th.\,9.7.1]{KS1}, $u$ may be represented by a single complex $F\in \rD^\rb_{\rc}(\pOS)$. We may choose a $\mu$-stratification $X=\bigsqcup_{\alpha}Z_{\alpha}$ such that $\shh^jF|_{Z_{\alpha}\times S}$ is isomorphic to $p^{-1}_{Z_{\alpha}}G_{\alpha, j}$, for some coherent $\sho_S$-module $G_{\alpha, j}$.

Let $X_k$ denote the (disjoint) union of the $k$-codimensional strata. By the same argument of
\cite[Th.\,9.7.1]{KS1}
we have $u=[F]=\sum_k[F_{X_k\times S}]$ and so $$[F]=\sum_{j,k}(-1)^j[(\shh^jF)_{X_k\times S}]$$
By assumption $\chi^S(F)=0$ hence
$$[\oplus_{j\, even} (\shh^j F)_{\{x\}\times S}]=[\oplus_{j\, odd} (\shh^jF)_{\{x\}\times S}]
$$
in $\rK_0(S)$, which implies
by
Remark~\ref{R1} 
$$[\oplus_{j\, even} (\shh^j F)_{Z_{\alpha}\times S}]=[\oplus_{j\, odd} (\shh^jF)_{Z_{\alpha}\times S}]$$
for each $\alpha$.
Because the $X_k$ are disjoint unions of $Z_{\alpha}$'s, we conclude that, for each $k$, 
$$[\oplus_{j\, even} (\shh^jF)_{X_k\times S}]=[\oplus_{j\, odd} (\shh^jF)_{X_k\times S}]$$ which implies that $[F]=0$.
\end{proof}

\begin{example}\label{exP}
Let us recall that given $S$ a smooth projective algebraic curve we have
$\rK_0(S)\simeq \mathbb Z\times \Pic(S)\simeq \mathbb Z^2\times \Pic_0(S)$
where $\Pic_0(S)$ is the kernel of the degree map. 
In particular
$\rK_0(\mathbb P^1(\C))$ is thus isomorphic to the unital ring $\Z^2$ with the usual addition and the multiplication  $(a,b)\cdot (c,d):= (ac, b+d)$ where, as expected, the unit is $(1,0)$ which corresponds to $[\sho_S]$. In particular, if to $F\in\rD^\rb_{\rc}(\pOS)$ corresponds $(a, b)\in\Z^2$ by the previous isomorphism $\chi^S$, to $[\bD F]$ corresponds $(a,-b)$.

 The datum of a $\mathbb P^1(\C)$-constructible function $f:X\to \rK_0(\mathbb P^1(\C))$ is equivalent to the datum of a locally finite sum 
$$\sum_{i,k}(z_i, z_k){\bf{I}}_{Z_{i,k}}$$ where each $(z_i, z_k)\in \Z^2$ and each $Z_{i,k}$ is compact subanalytic contractible in $X$.
\end{example}

\providecommand{\sortnoop}[1]{}\providecommand{\eprint}[1]{\href{http://arxiv.org/abs/#1}{\texttt{arXiv\string:\allowbreak#1}}}\providecommand{\hal}[1]{\href{https://hal.archives-ouvertes.fr/hal-#1}{\texttt{hal-#1}}}
\subsection*{Acknowledgements}
{The research of L.~Fiorot was supported by Project funded by the EuropeanUnion – NextGenerationEU under the National Recovery and Resilience Plan (NRRP), Mission 4 Component 2 Investment 1.1 - Call PRIN 2022 No. 104 of February 2, 2022 of Italian Ministry of University and Research; Project 2022S97PMY (subject area: PE - Physical Sciences and Engineering) ``Structures for Quivers, Algebras and Representations (SQUARE)''. L. Fiorot is moreover member of INDAM - GNSAGA.
The research of T.~Monteiro Fernandes was supported by
Funda\c c{\~a}o para a Ci{\^e}ncia e Tecnologia, under the project: UIDB/04561/2020}.

\end{document}